\documentclass[psamsfonts]{amsart}
\usepackage{amsmath}
\usepackage{amssymb}
\usepackage{mathtools}
\usepackage{pifont}
\usepackage{multicol}
\usepackage{mathrsfs}
\usepackage{multicol}
\usepackage{amsthm}
\usepackage{float}
\usepackage{verbatim} 
\usepackage{tikz-cd}
\usetikzlibrary{arrows}
\usetikzlibrary{shapes.arrows}   
\usetikzlibrary{positioning}
\usepackage{mathrsfs}
\usepackage{lscape} 
\usepackage{wasysym}
\usepackage{mathtools}
\usepackage{python}
\usepackage{enumerate}
\usepackage[foot]{amsaddr}
\usepackage[numbers,comma,sort&compress]{natbib} 
\usepackage[colorlinks=true,citecolor=blue]{hyperref}

\usepackage{adjustbox}
\usepackage{tikz}
\usepackage{tikz-3dplot}
\usepackage{makecell}
\usepackage{caption}
\usetikzlibrary{shapes,calc,positioning,intersections}
\tdplotsetmaincoords{80}{120}

\definecolor{1}{rgb}{1,0.2,0.3}
\definecolor{2}{rgb}{0.1,0.3,0.5}
\definecolor{3}{rgb}{1,1,0}
\definecolor{4}{rgb}{255,255,255}
\newcommand{\R}{\mathbb{R}}	   

\newcommand{\Z}{\mathbb{Z}}		   

\usepackage{mathtools}

\def\co{\colon\thinspace}

\usepackage{lipsum}

\newtheorem{theorem}{Theorem}[section]
\newtheorem{prop}[theorem]{Proposition}
\newtheorem{corollary}[theorem]{Corollary}
\newtheorem{lemma}[theorem]{Lemma}
\theoremstyle{definition}
\newtheorem{definition}[theorem]{Definition}
\newtheorem{remark}[theorem]{Remark}
\newtheorem{example}[theorem]{Example}

\newtheorem{conj}[theorem]{Conjecture}

\begin{document}
	\title{Higher-dimensional cubical sliding puzzles}

	\author{Moritz Beyer}
	\email{moritzbeyer@t-online.de}

	\author{Stefano Mereta}
	\email{stefano.mereta@mis.mpg.de}

	\author{Érika Roldán}
	\email{erika.roldan@mis.mpg.de} 
	\address{Max Planck Institute for Mathematics in the Sciences, Inselstraße 22, 04103 Leipzig, Germany}
	\address{ScaDS.AI Leipzig, Humboldtstraße 25, 04105 Leipzig}

	\author{Peter Voran}
	\email{p-voran@t-online.de}

	\date{\today}
	

	\maketitle 
	\begin{abstract}
	We introduce higher-dimensional cubical sliding puzzles that are inspired by the classical 15 Puzzle from the 1880s. In our puzzles, on a $d$-dimensional cube, a labeled token can be slid from one vertex to another if it is \textit{topologically free to move} on lower-dimensional faces. We analyze solvability of these puzzles by studying how the puzzle graph changes with the number of labeled tokens vs empty vertices. We give characterizations of the different regimes ranging from being completely stuck (and thus all puzzles unsolvable) to having only one giant component where almost all puzzles can be solved. For the Cube, the Tesseract, and the Penteract ($5$-dimensional cube) we have implemented an algorithm to completely analyze their solvability and we provide specific puzzles for which we even know the minimum number of moves needed to solve them. Readers are encouraged to try to solve some of our puzzles on this \href{https://www.erikaroldan.net/cubical-sliding-puzzles}{website}  \cite{videogame}.  \end{abstract}
 
	\section{Introduction}\label{sec:introduction}

One of the most famous puzzles of all time is the 15 puzzle which consists of a $4 \times 4$ square board divided into 16 unit squares with 15 labeled tiles from 1 to 15 placed on the board and with the one remaining square being empty. The initial configuration of the tiles in the 15 puzzle has the tiles 1 to 13 placed in ascending order (from left to right and top to bottom), and the tiles 14 and 15 interchanged. The target configuration of this puzzle leaves the first 13 ordered tiles fixed and interchanges the 14 and 15 tiles. In 1879 \cite{johnson1879notes}, Johnson proved that this puzzle has no solution. In the second part of the same paper \cite{story1879notes}, Story further proved that the set of all possible configurations of the labeled tiles on the board is divided into two giant components which determine which puzzles are solvable or not \footnote{By a puzzle we mean a fixed initial and a fixed target configuration.}.

 Other sliding puzzles based on regular tessellations include those that have more general polyominoes as the pieces to be slid \cite{hordern}, or with boards consisting of a subset of tiles from the regular hexagonal tessellation of the plane with hexagonal labeled pieces \cite{alpert2020discrete, hex}. 

Here, we introduce a new family of sliding puzzles whose boards are the graphs of the $d$-cube with labeled tokens placed at their vertices as the pieces to be slid under certain topological obstructions. Of course, if all the cube's vertices are occupied by labeled tokens then the tokens cannot move at all. In this case, any sliding puzzle is clearly unsolvable. If some tokens of the board are removed, then empty space is created, and some sliding puzzles may become solvable. 

Before formally stating the sliding rules, we encourage the curious reader to try to solve some of these cubical sliding puzzles on this \href{https://www.erikaroldan.net/cubical-sliding-puzzles}{website} \cite{videogame}. We consider deducing the topological obstructions that govern the movement of the tokens a nice first riddle to solve. 

Intuitively, for each fixed dimension $k\leq d$, a token on a vertex of a $k$-dimensional face of the $d$-cube requires all other vertices of this $k$-face to be empty to be able to slide or move to any other unoccupied vertex of the same $k$-face. For instance, if $k=1$, a token can  move to any empty vertex that shares an edge with the vertex where the token is placed. This $k=1$ puzzle model, which in particular applies to the 15 Puzzle, has been previously studied for more general graphs, and the solvability of these puzzles is very much understood \cite{wilson}.

To give a more precise definition of the puzzles and the rules for movement, we first need to set some notation. All along the paper, $k,d$ and $l$ will denote integer numbers with the following constraints: $d \ge 1$, $1 \le l \le 2^d$, and $k \le d$. In general, $d$ denotes the dimension of the cube on which we build a sliding puzzle with $2^d-l$ labeled tokens placed at its vertices and the remaining $l$ vertices being unoccupied. The number $k$ determines the rule for moving a token to an unoccupied vertex. For a given $d$, we denote by $Q^d$ the graph containing all the vertices and edges of the $d$-dimensional cube.

Let $\mathcal{C}[d,l]$ be the set of all configurations of $2^d-l$ labeled tokens on the vertices of $Q^d$, and $\mathcal{L}_C =\{1,2,...,2^d-l\}$ the set of labels of these tokens.

\begin{definition}
	Let $\sigma$ be a $k$-dimensional face of $Q^d$ and $t$ a labeled token on one of the vertices of $\sigma$. If the remaining vertices of $\sigma$ are unoccupied, then we say that $t$ is in a \textit{free $k$-state with respect to $\sigma$}. A labeled token that is on a free $k$-state with respect to  $\sigma$ can move (or slide) to occupy any of the vertices of $\sigma$. We call this move of $t$ a \textit{$k$-move on $\sigma$}. 
	
\end{definition}

\begin{definition}
	For a triple $(d,k,l)$, let the $(d,k,l)$-th puzzle graph, denoted as $puz[d,k,l]$, be the graph whose set of vertices is $\mathcal{C}[d,l]$ and whose set of edges is given by
	\[
	\{(C,D) \in \mathcal{C}[d,l] \times \mathcal{C}[d,l] \mid \text {$C$ can be reached from $D$ by a single $k$-move} \}.
	\]
\end{definition}	
Observe that this is a symmetric relation, that is, $C$ can be reached from $D$ by a single $k$-move if and only if $D$ can be reached from $C$ by a single $k$-move.

\begin{definition}\label{def:semi-isolated}
    We say that a configuration $C$ is:
    \begin{enumerate}
        \item $k$\emph{-isolated}, if no token can perform a $k$-move;
        \item $k$\emph{-semi-isolated}, if there is a non-empty subset $\mathcal A_C \subsetneq L_C$ such that each labeled configuration in the connected component of $C$ in $puz[d,k,l]$ is obtained by a series of $k$-moves affecting only elements of $\mathcal L_C \setminus \mathcal A_C$. We refer to this set $\mathcal A_C$ as the \emph{set of stuck tokens} of $C$ and to any token in $\mathcal A_C$ as a \emph{stuck token} of $C$;
        \item $k$\emph{-mobile}, if it is neither $k$-isolated or $k$-semi-isolated.
    \end{enumerate}
\end{definition}
\begin{example}
The configuration in the middle in Figure \ref{fig:example-2-semi-isolated} is a $2$-semi-isolated configuration. With this $k=2$ rule, only the token labeled with 5 can move. In this case, the set of stuck tokens is $ \mathcal A_C=\{1, 2, 3, 4\}$.
\end{example}
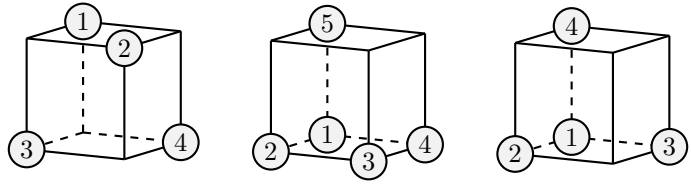
\begin{figure}
    \centering
\begin{tikzpicture}[scale=1.5, tdplot_main_coords,axis/.style={->},thick]

  \draw[dashed,opacity=1] (0,0,0) -- (0,1,0);
      \draw[dashed,opacity=1] (0,0,0) -- (1,0,0);
      \draw[dashed,opacity=1] (0,0,0) -- (0,0,1);
      \draw[thick,opacity=1] (0,1,0) -- (0,1,1);
      \draw[thick,opacity=1] (0,1,0) -- (1,1,0);
      \draw[thick,opacity=1] (1,0,0) -- (1,1,0);
      \draw[thick,opacity=1] (1,0,0) -- (1,0,1);
      \draw[thick,opacity=1] (0,0,1) -- (0,1,1);
      \draw[thick,opacity=1] (0,0,1) -- (1,0,1); 
      \draw[thick,opacity=1] (1,1,0) -- (1,1,1);
      \draw[thick,opacity=1] (0,1,1) -- (1,1,1);
      \draw[thick,opacity=1] (1,0,1) -- (1,1,1);

\draw[fill=gray!10] (1,0,0) circle (0.45em); 
\draw[fill=gray!10] (0,1,0) circle (0.45em); 
\draw[fill=gray!10] (0,0,1) circle (0.45em); 
\draw[fill=gray!10] (1,1,1) circle (0.45em); 

\node at (0,1,0) {$4$};
\node at (1,0,0) {$3$};
\node at (0,0,1) {$1$};
\node at (1,1,1) {$2$};

\tikzset{shift={(0,2.5,0.2)}}

\draw[dashed,opacity=1] (0,0,0) -- (0,1,0);
      \draw[dashed,opacity=1] (0,0,0) -- (1,0,0);
      \draw[dashed,opacity=1] (0,0,0) -- (0,0,1);
      \draw[thick,opacity=1] (0,1,0) -- (0,1,1);
      \draw[thick,opacity=1] (0,1,0) -- (1,1,0);
      \draw[thick,opacity=1] (1,0,0) -- (1,1,0);
      \draw[thick,opacity=1] (1,0,0) -- (1,0,1);
      \draw[thick,opacity=1] (0,0,1) -- (0,1,1);
      \draw[thick,opacity=1] (0,0,1) -- (1,0,1); 
      \draw[thick,opacity=1] (1,1,0) -- (1,1,1);
      \draw[thick,opacity=1] (0,1,1) -- (1,1,1);
      \draw[thick,opacity=1] (1,0,1) -- (1,1,1);

\draw[fill=gray!10] (0,0,0) circle (0.45em); 
\draw[fill=gray!10] (1,0,0) circle (0.45em); 
\draw[fill=gray!10] (0,1,0) circle (0.45em); 
\draw[fill=gray!10] (0,0,1) circle (0.45em); 
\draw[fill=gray!10] (1,1,0) circle (0.45em); 

\node at (0,0,0) {$1$};
\node at (0,1,0) {$4$};
\node at (1,0,0) {$2$};
\node at (0,0,1) {$5$};
\node at (1,1,0) {$3$};

\tikzset{shift={(0,2.5,0.2)}}

 \draw[dashed,opacity=1] (0,0,0) -- (0,1,0);
      \draw[dashed,opacity=1] (0,0,0) -- (1,0,0);
      \draw[dashed,opacity=1] (0,0,0) -- (0,0,1);
      \draw[thick,opacity=1] (0,1,0) -- (0,1,1);
      \draw[thick,opacity=1] (0,1,0) -- (1,1,0);
      \draw[thick,opacity=1] (1,0,0) -- (1,1,0);
      \draw[thick,opacity=1] (1,0,0) -- (1,0,1);
      \draw[thick,opacity=1] (0,0,1) -- (0,1,1);
      \draw[thick,opacity=1] (0,0,1) -- (1,0,1); 
      \draw[thick,opacity=1] (1,1,0) -- (1,1,1);
      \draw[thick,opacity=1] (0,1,1) -- (1,1,1);
      \draw[thick,opacity=1] (1,0,1) -- (1,1,1);

\draw[fill=gray!10] (0,0,0) circle (0.45em); 
\draw[fill=gray!10] (1,0,0) circle (0.45em); 
\draw[fill=gray!10] (0,0,1) circle (0.45em); 
\draw[fill=gray!10] (0,1,0) circle (0.45em); 

\node at (0,0,0) {$1$};
\node at (0,1,0) {$3$};
\node at (1,0,0) {$2$};
\node at (0,0,1) {$4$};

\end{tikzpicture}
    \caption{From left to right: a $2$-isolated, a $2$-semi-isolated and a $2$-mobile configuration on the Cube.}
    \label{fig:example-2-semi-isolated}
\end{figure}

\begin{definition}
    Let $s(d,k)$ be the minimum natural number such that $puz[d,k,s(d,k)]$ contains a configuration that is not $k$-isolated i.e.\ such that there exists at least one $k$-semi-isolated or one $k$-mobile configuration. 
\end{definition}

	 A $k$ dimensional face $\sigma$ in $Q^d$ has $2^k$ vertices. Thus, by definition $s(d,k)=2^k-1$. The value $s(d,k)$ captures when there is enough space, i.e.\ vertices without labeled tokens, for the tokens to move. A natural question to ask is if there is a regime in which there exist only configurations that are either $k$-isolated or $k$-semi-isolated (that is, without $k$-mobile configurations). In the case of $k = 1$ and $l=1$, it was proven by Wilson \cite{wilson} that $puz[d,1,s(d,k)]$ contains no $k$-isolated or $k$-semi-isolated but only $1$-mobile configurations, moreover, the author proved that $puz[d,1,s(d,k)]$ has exactly two isomorphic connected components. In our first main result, we prove that this is not always the case for puzzle graphs of cubical sliding puzzles if $k\geq 2$. In other words, there is a regime in which for some values of $l$, that depend on $k$ and $d$, the puzzle graph has at least one $k$-semi-isolated configuration and no components consisting of $k$-mobile configurations.  The precise statement is stated below in Theorem \ref{thm:k-semi-isolated regime}.

\begin{theorem} \label{thm:k-semi-isolated regime}
If $d > k \geq 2$, the graph $puz[d,k,s(d,k)]$  contains a $k$-semi-isolated configuration and no $k$-mobile configurations.
\end{theorem}

    Through the paper, we say that a configuration $C$ is obtained by performing a permutation $g$ of another configuration $D$, if the underlying unlabeled configuration is the same for $C$ and $D$, and the labels are permuted by $g$.
  
\begin{definition}
	We say that $puz[d, k, l]$  has the \emph{even solvability property} if there exists a  configuration $C$, such that all configurations that can be obtained from $C$ by an even permutation are in the same connected component of $puz[d,k,l]$ as $C$. In this case, we will say that such a configuration $C$ has the even solvability property.
\end{definition}
We will prove in Lemma \ref{lemma:even-solvability-k-mobile} that if a configuration  has the even solvability property, then that configuration is $k$-mobile.

\begin{definition}
    For any pair of integers $(d,k)$ let $S(d,k)$ be the minimum integer such that $puz[d,k,S(d,k)]$ has the even solvability property.
\end{definition}

Our second main result is to prove that while $l$ increases, the even solvability property is achieved exactly when the first $k$-mobile configuration appears.

\begin{theorem}\label{thm:S(d,k) = smallest-containing-not-k-semi-isolated}
 For any admissible tuple of integers $(d, k)$, the value $S(d,k)$ is the smallest number of unoccupied vertices such that $puz[d, k, S(d,k)]$ contains a $k$-mobile configuration. Moreover, the puzzle graph $puz[ d , k , S(d,k) ]$ has either one or two connected components containing $k$-mobile configurations.     
\end{theorem}

In the case of two connected components containing $k$-mobile configurations, those are separated only by parity: for each configuration $C$ all its even permutations are contained in one of them, all the odd permutations in the other. This phenomenon of having \textit{two giant components} determined by parity is what governs solvability in the 15 Puzzle. It was recently proven that in hexagonal sliding puzzles it can happen that the stage of having two giant components is not always attained\cite{hex}. Following the convention established in \cite{hex}, we say that a puzzle graph has the \textit{strong parity property} if it has exactly two components containing $k$-mobile configurations.

From Theorem \ref{thm:S(d,k) = smallest-containing-not-k-semi-isolated} above and Lemma \ref{lemma:even-solvability-k-mobile}, we get the following Corollary:
\begin{corollary}
 For any admissible tuple of integers $(d, k)$, if $l < S(d,k)$ the puzzle graph $puz[d,k,l]$ contains only $k$-isolated or $k$-semi-isolated configurations.   
\end{corollary} 

If we increase the number of unoccupied vertices beyond the critical point $S(d,k)$ the puzzle graph has only one component containing $k$-mobile configurations. We prove this in Theorem \ref{theorem:maximality-of-L-introduction} stated below. 

\begin{theorem}\label{theorem:maximality-of-L-introduction}
	For any admissible tuple of integers $(d, k)$, the subgraph of the puzzle graph $puz[d,k,S(d,k)+1]$ whose vertices are all the $k$-mobile configurations is connected.
\end{theorem}

Motivated by Theorem \ref{theorem:maximality-of-L-introduction} and the computations that we present in Section \ref{sec:computational-results}, we conjecture that:

\begin{conj}\label{conjecture:strong-parity}
     The puzzle graph $puz[d,k,S(d,k)]$ has the strong parity property. That is, it has exactly two components containing $k$-mobile configurations.
\end{conj}

So far, our results aim to give a fairly complete picture of how the structure and the dynamics of the puzzle graph $puz[d,k,l]$ change accordingly with $l$, the number of unoccupied vertices of a cubical sliding puzzle. In the following two results, we give exact expressions for $S(d,k)$.

\begin{theorem}\label{theorem:value-of-L(k+1,k)}
	For any $k \ge 1$,
	\[
	S(k+1,k) = \sum_{i = 1}^{k} (2^i - 1) = 2^{k+1} - k - 2
	\]
i.e.\ the even solvability property is achieved with $k+2$ labeled tokens.
\end{theorem}

Starting from this first formula, in Section \ref{sec:values-of-S(d,k)} we prove a recursive and exact formula for $S(d,k)$ for all admissible pairs $(d,k)$.

\begin{theorem}\label{theorem:d=k+n-intro}
For any $k \ge 2$ and $n \ge 2$ we have:
\[
S(n+k , k) = S(n+k - 1, k-1) + S(n+k - 1, k).
\]
this expands to:
\[
S(n+k , k) = \sum_{i = 2}^{k} \binom{n +k-i}{k-i} (2^{i} - 1) + \sum_{i = 2}^{n} \binom{n+k-i}{k-1}.
\]
\end{theorem}

Finally, we analyze in detail the puzzle graphs for cubical sliding puzzles defined in the Cube, the Tesseract, and the Penteract. In particular, we are interested in studying the diameter of these puzzle graphs. To motivate our results, we first give an example of an optimal solution to a sliding puzzle on the Cube with the $k = 2$ rule.

 \begin{example}
    For $d = 3$ and $k = 2$ the solution to a sliding puzzle is given by the following sequence of moves (as we will discuss later, by applying a $3$-cycle on a configuration): 
    
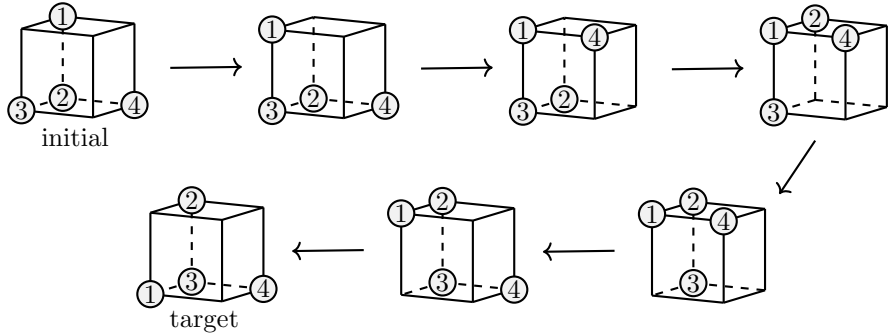
\vspace{0.3cm}

    \begin{tikzpicture}[scale=1.1, tdplot_main_coords,axis/.style={->},thick]

      \draw[dashed,opacity=1] (0,0,0) -- (0,1,0);
      \draw[dashed,opacity=1] (0,0,0) -- (1,0,0);
      \draw[dashed,opacity=1] (0,0,0) -- (0,0,1);
      \draw[thick,opacity=1] (0,1,0) -- (0,1,1);
      \draw[thick,opacity=1] (0,1,0) -- (1,1,0);
      \draw[thick,opacity=1] (1,0,0) -- (1,1,0);
      \draw[thick,opacity=1] (1,0,0) -- (1,0,1);
      \draw[thick,opacity=1] (0,0,1) -- (0,1,1);
      \draw[thick,opacity=1] (0,0,1) -- (1,0,1); 
      \draw[thick,opacity=1] (1,1,0) -- (1,1,1);
      \draw[thick,opacity=1] (0,1,1) -- (1,1,1);
      \draw[thick,opacity=1] (1,0,1) -- (1,1,1);

\draw[fill=gray!10] (0,0,0) circle (0.45em); 
\draw[fill=gray!10] (1,0,0) circle (0.45em); 
\draw[fill=gray!10] (0,1,0) circle (0.45em); 
\draw[fill=gray!10] (0,0,1) circle (0.45em); 

\node at (0,0,0) {$2$};
\node at (0,1,0) {$4$};
\node at (1,0,0) {$3$};
\node at (0,0,1) {$1$};
\node at (2.3,1.5,0) {initial};

\draw[->,opacity=1] (0,1.5,0.5) -- (0,2.5,0.6);

\tikzset{shift={(0,3.5,0.3)}}

  \draw[dashed,opacity=1] (0,0,0) -- (0,1,0);
      \draw[dashed,opacity=1] (0,0,0) -- (1,0,0);
      \draw[dashed,opacity=1] (0,0,0) -- (0,0,1);
      \draw[thick,opacity=1] (0,1,0) -- (0,1,1);
      \draw[thick,opacity=1] (0,1,0) -- (1,1,0);
      \draw[thick,opacity=1] (1,0,0) -- (1,1,0);
      \draw[thick,opacity=1] (1,0,0) -- (1,0,1);
      \draw[thick,opacity=1] (0,0,1) -- (0,1,1);
      \draw[thick,opacity=1] (0,0,1) -- (1,0,1); 
      \draw[thick,opacity=1] (1,1,0) -- (1,1,1);
      \draw[thick,opacity=1] (0,1,1) -- (1,1,1);
      \draw[thick,opacity=1] (1,0,1) -- (1,1,1);

\draw[fill=gray!10] (0,0,0) circle (0.45em); 
\draw[fill=gray!10] (1,0,0) circle (0.45em); 
\draw[fill=gray!10] (0,1,0) circle (0.45em); 
\draw[fill=gray!10] (1,0,1) circle (0.45em); 

\node at (0,0,0) {$2$};
\node at (0,1,0) {$4$};
\node at (1,0,0) {$3$};
\node at (1,0,1) {$1$};

\draw[ ->,opacity=1] (0,1.5,0.5) -- (0,2.5,0.6);

\tikzset{shift={(0,3.5,0.3)}}

  \draw[dashed,opacity=1] (0,0,0) -- (0,1,0);
      \draw[dashed,opacity=1] (0,0,0) -- (1,0,0);
      \draw[dashed,opacity=1] (0,0,0) -- (0,0,1);
      \draw[thick,opacity=1] (0,1,0) -- (0,1,1);
      \draw[thick,opacity=1] (0,1,0) -- (1,1,0);
      \draw[thick,opacity=1] (1,0,0) -- (1,1,0);
      \draw[thick,opacity=1] (1,0,0) -- (1,0,1);
      \draw[thick,opacity=1] (0,0,1) -- (0,1,1);
      \draw[thick,opacity=1] (0,0,1) -- (1,0,1); 
      \draw[thick,opacity=1] (1,1,0) -- (1,1,1);
      \draw[thick,opacity=1] (0,1,1) -- (1,1,1);
      \draw[thick,opacity=1] (1,0,1) -- (1,1,1);

\draw[fill=gray!10] (0,0,0) circle (0.45em); 
\draw[fill=gray!10] (1,0,0) circle (0.45em); 
\draw[fill=gray!10] (1,1,1) circle (0.45em); 
\draw[fill=gray!10] (1,0,1) circle (0.45em); 

\node at (0,0,0) {$2$};
\node at (1,1,1) {$4$};
\node at (1,0,0) {$3$};
\node at (1,0,1) {$1$};

\draw[ ->,opacity=1] (0,1.5,0.5) -- (0,2.5,0.6);

\tikzset{shift={(0,3.5,0.3)}}

  \draw[dashed,opacity=1] (0,0,0) -- (0,1,0);
      \draw[dashed,opacity=1] (0,0,0) -- (1,0,0);
      \draw[dashed,opacity=1] (0,0,0) -- (0,0,1);
      \draw[thick,opacity=1] (0,1,0) -- (0,1,1);
      \draw[thick,opacity=1] (0,1,0) -- (1,1,0);
      \draw[thick,opacity=1] (1,0,0) -- (1,1,0);
      \draw[thick,opacity=1] (1,0,0) -- (1,0,1);
      \draw[thick,opacity=1] (0,0,1) -- (0,1,1);
      \draw[thick,opacity=1] (0,0,1) -- (1,0,1); 
      \draw[thick,opacity=1] (1,1,0) -- (1,1,1);
      \draw[thick,opacity=1] (0,1,1) -- (1,1,1);
      \draw[thick,opacity=1] (1,0,1) -- (1,1,1);

\draw[fill=gray!10] (0,0,1) circle (0.45em); 
\draw[fill=gray!10] (1,0,0) circle (0.45em); 
\draw[fill=gray!10] (1,1,1) circle (0.45em); 
\draw[fill=gray!10] (1,0,1) circle (0.45em); 

\node at (0,0,1) {$2$};
\node at (1,1,1) {$4$};
\node at (1,0,0) {$3$};
\node at (1,0,1) {$1$};

\draw[->,opacity=1] (0,0,-0.5) -- (0,-0.5,-1.2);

\tikzset{shift={(0,-1.7,-2.4)}}

  \draw[dashed,opacity=1] (0,0,0) -- (0,1,0);
      \draw[dashed,opacity=1] (0,0,0) -- (1,0,0);
      \draw[dashed,opacity=1] (0,0,0) -- (0,0,1);
      \draw[thick,opacity=1] (0,1,0) -- (0,1,1);
      \draw[thick,opacity=1] (0,1,0) -- (1,1,0);
      \draw[thick,opacity=1] (1,0,0) -- (1,1,0);
      \draw[thick,opacity=1] (1,0,0) -- (1,0,1);
      \draw[thick,opacity=1] (0,0,1) -- (0,1,1);
      \draw[thick,opacity=1] (0,0,1) -- (1,0,1); 
      \draw[thick,opacity=1] (1,1,0) -- (1,1,1);
      \draw[thick,opacity=1] (0,1,1) -- (1,1,1);
      \draw[thick,opacity=1] (1,0,1) -- (1,1,1);

\draw[fill=gray!10] (0,0,1) circle (0.45em); 
\draw[fill=gray!10] (0,0,0) circle (0.45em); 
\draw[fill=gray!10] (1,1,1) circle (0.45em); 
\draw[fill=gray!10] (1,0,1) circle (0.45em); 

\node at (0,0,1) {$2$};
\node at (1,1,1) {$4$};
\node at (0,0,0) {$3$};
\node at (1,0,1) {$1$};

\draw[ ->,opacity=1] (0,-1.1,0.3) -- (0,-2.1,0.2);

\tikzset{shift={(0,-3.5,-0.3)}}

  \draw[dashed,opacity=1] (0,0,0) -- (0,1,0);
      \draw[dashed,opacity=1] (0,0,0) -- (1,0,0);
      \draw[dashed,opacity=1] (0,0,0) -- (0,0,1);
      \draw[thick,opacity=1] (0,1,0) -- (0,1,1);
      \draw[thick,opacity=1] (0,1,0) -- (1,1,0);
      \draw[thick,opacity=1] (1,0,0) -- (1,1,0);
      \draw[thick,opacity=1] (1,0,0) -- (1,0,1);
      \draw[thick,opacity=1] (0,0,1) -- (0,1,1);
      \draw[thick,opacity=1] (0,0,1) -- (1,0,1); 
      \draw[thick,opacity=1] (1,1,0) -- (1,1,1);
      \draw[thick,opacity=1] (0,1,1) -- (1,1,1);
      \draw[thick,opacity=1] (1,0,1) -- (1,1,1);

\draw[fill=gray!10] (0,0,1) circle (0.45em); 
\draw[fill=gray!10] (0,0,0) circle (0.45em); 
\draw[fill=gray!10] (0,1,0) circle (0.45em); 
\draw[fill=gray!10] (1,0,1) circle (0.45em); 

\node at (0,0,1) {$2$};
\node at (0,1,0) {$4$};
\node at (0,0,0) {$3$};
\node at (1,0,1) {$1$};

\draw[ ->,opacity=1] (0,-1.1,0.3) -- (0,-2.1,0.2);

\tikzset{shift={(0,-3.5,-0.3)}}

  \draw[dashed,opacity=1] (0,0,0) -- (0,1,0);
      \draw[dashed,opacity=1] (0,0,0) -- (1,0,0);
      \draw[dashed,opacity=1] (0,0,0) -- (0,0,1);
      \draw[thick,opacity=1] (0,1,0) -- (0,1,1);
      \draw[thick,opacity=1] (0,1,0) -- (1,1,0);
      \draw[thick,opacity=1] (1,0,0) -- (1,1,0);
      \draw[thick,opacity=1] (1,0,0) -- (1,0,1);
      \draw[thick,opacity=1] (0,0,1) -- (0,1,1);
      \draw[thick,opacity=1] (0,0,1) -- (1,0,1); 
      \draw[thick,opacity=1] (1,1,0) -- (1,1,1);
      \draw[thick,opacity=1] (0,1,1) -- (1,1,1);
      \draw[thick,opacity=1] (1,0,1) -- (1,1,1);

\draw[fill=gray!10] (0,0,1) circle (0.45em); 
\draw[fill=gray!10] (0,0,0) circle (0.45em); 
\draw[fill=gray!10] (0,1,0) circle (0.45em); 
\draw[fill=gray!10] (1,0,0) circle (0.45em); 

\node at (0,0,1) {$2$};
\node at (0,1,0) {$4$};
\node at (0,0,0) {$3$};
\node at (1,0,0) {$1$};
\node at (2.3,1.5,0) {target};

\end{tikzpicture}

\vspace{0.3cm}

    This sequence of movements corresponds to a unique shortest path in the puzzle graph.
    To see why there cannot exist any shorter path, note that each of the tokens 1, 2, and 3 have to move at least once to reach their desired location. Also, one of them has to move at least once more, because at the start all desired locations are occupied. Moreover, it is easy to see that it is not possible to solve the puzzle without moving in the last move the token labeled with 4, which thus has to move at least twice.  It follows that the puzzle cannot be solved in less than 6 moves.
\end{example}

The following results are based on the computations of the puzzle graphs for almost all values of $l$ of the Cube, the Tesseract, and some for the Penteract that we present in Section \ref{sec:computational-results}. 

\begin{theorem}\label{small-dimension-3}
    For $d = 3$ and $k = 2$ we have $S(3,2) = 4$ and any puzzle can be solved in at most $10$ moves. Furthermore strong parity property holds for $puz[3,2,S(3,2)]$.
 \end{theorem}

\begin{theorem}\label{small-dimension-4}
    For $d = 4$ and $k = 3$ we have $S(4,3) = 11$ and any puzzle can be solved in at most $17$ moves. Furthermore strong parity property holds for $puz[4,3,S(4,3)]$.
 \end{theorem}

\begin{theorem}\label{small-dimension-4-2}
      For $d = 4$ and $k = 2$ we have $S(4,2) = 5$. Again strong parity property holds for $puz[4,2,S(4,2)]$.
\end{theorem}

It is known that computing the diameter of sliding puzzles with the $1$-rule on general graphs is NP-hard \cite{goldreich}. However, based on the computations that we present in Section \ref{sec:computational-results}, the following pattern seems to hold:

\begin{conj}
    For any $l \in \{2,3,\dots,2^d\}$ the diameter of $puz[d,1,l]$ is given by $diam=d \cdot (2^d-l)$.
\end{conj}

\subsection{Structure of the paper}
In Section \ref{sec:preliminary} we give preliminary results and we prove Theorem \ref{thm:k-semi-isolated regime}. Section \ref{sec:unlabeled puzzle graph} contains some technical results about the unable puzzle graph that we require for some proofs and our computational results. 
In Section \ref{sec:values-of-S(d,k)}, we present a proof of Theorem \ref{theorem:value-of-L(k+1,k)} and of Theorem \ref{theorem:d=k+n-intro}. In Section \ref{sec:k-mobile-configurations}, we prove 
Theorem \ref{thm:S(d,k) = smallest-containing-not-k-semi-isolated} and Theorem \ref{theorem:maximality-of-L-introduction}.
Lastly, in Section \ref{sec:computational-results}, we present the algorithms used for the computations, and using these results we get Theorems \ref{small-dimension-3}, \ref{small-dimension-4}, and \ref{small-dimension-4-2}.

\section{Definitions and first results}\label{sec:preliminary}

We start this section with a first observation about the structure of puzzle graphs giving a lower bound on the degree of non-isolated vertices.

\begin{lemma}\label{degree of vertices}
Every vertex $v \in puz[d,k,l]$ corresponding to a non-isolated configuration has degree at least $2^k-1$.
\end{lemma}
\begin{proof}
The vertex $v \in puz[d,k,l]$ corresponds to a configuration that has at least a token that can perform a $k$-move. This token can perform at least $2^k-1$ different moves. This results in $2^k-1$ different configurations which are connected to $v$. Thus $deg(v) \geq 2^k-1$.
\end{proof}

We now prove Theorem \ref{thm:k-semi-isolated regime} from the Introduction.

\begin{theorem}
If $d > k \geq 2$, $puz[d,k,s(d,k)]$  contains a $k$-semi-isolated configuration and no $k$-mobile configurations.
\end{theorem}

\begin{proof}
    By definition of $s(d,k)$, there exists a configuration a $k$-face $\sigma$ of which contains one token $t$ and $2^k - 1$ unoccupied vertices. This configuration is not $k$-isolated as $t$ can move. It is not $k$-mobile, because as $d > k$ there are $2^d - 2^k$ other tokens. Those tokens define a stuck set, as any $k$-face other than $\sigma$ contains no more than $2^{k-1}$ unoccupied vertices and thus no movement is allowed outside of $\sigma$.
\end{proof}
For $d = k$ the number $s(k,k)$ is equal to $2^k - 1$. Thus there is only one token. This token can move freely, every configuration is $k$-mobile and $puz(k, k, 2^k - 1)$ is connected.
\begin{lemma}\label{lemma:even-solvability-k-mobile}
Let $d,k,l$ such that $puz(d,k,l)$ has the even solvability property. Then $puz(d,k,l)$ contains a $k$-mobile configuration.
\end{lemma}
\begin{proof}
    By definition $puz(d,k,l)$ contains a configuration $C$  such that all configurations that can be obtained from $C$ by an even permutation are in the same connected component of $puz(d, k, l)$. This configuration $C$ is $k$-mobile, as otherwise there would exist a stuck set $A$ and any permutation of $C$ not fixing at least one token in $A$ cannot be obtained from $C$, as this token is stuck.
\end{proof}
\subsection{Lifts of configurations}\label{lifts} In this subsection we introduce the notion of lift of a configuration and discuss some of its properties.
\begin{definition}
        Given an integer $n \ge 0$, let 
		\[
			q_{d,l}(n) := l + \sum_{i=0}^{n-1} 2^{d+i} - n
		\]
Given a configuration $C \in \mathcal{C}[d,l]$, a labeled configuration $C' \in \mathcal{C}[d+1,q_{d,l}(1)]$ is a \textit{$1$-lift} of $C$, if it is obtained by choosing an embedding of $Q^d$ as a face of $Q^{d+1}$ and by positioning an additional labeled token on a vertex of the face parallel to the image of this embedding. Given any $n \geq 0$ and a configuration $C \in \mathcal{C}[d,l]$, we will say that a labeled configuration $C' \in \mathcal{C}[d+n, q_{d,l}(n)]$ is an \textit{$n$-lift} of $C$ if it can be obtained from $C$ by $n$ iterated lifts. 
\end{definition}

\begin{example}
The following is an example of a $1$-lift of a configuration from $\mathcal{C}[3,4]$ to $\mathcal{C}[4,11]$. 
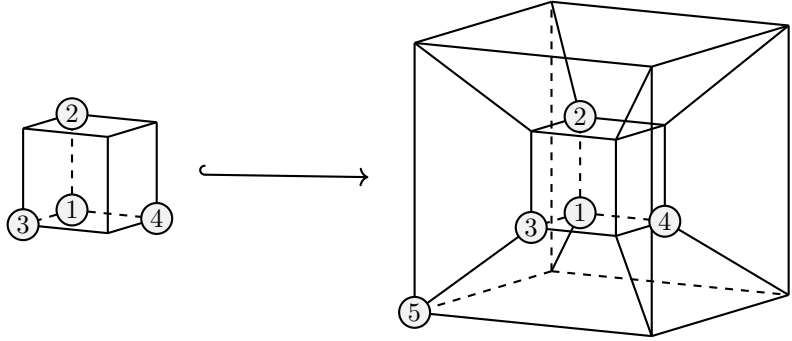
\vspace{0.25cm}

\begin{center}
\begin{tikzpicture}[scale=1.3, tdplot_main_coords,axis/.style={->},thick]

      \draw[dashed,opacity=1] (0,0,0) -- (0,1,0);
      \draw[dashed,opacity=1] (0,0,0) -- (1,0,0);
      \draw[dashed,opacity=1] (0,0,0) -- (0,0,1);
      \draw[thick,opacity=1] (0,1,0) -- (0,1,1);
      \draw[thick,opacity=1] (0,1,0) -- (1,1,0);
      \draw[thick,opacity=1] (1,0,0) -- (1,1,0);
      \draw[thick,opacity=1] (1,0,0) -- (1,0,1);
      \draw[thick,opacity=1] (0,0,1) -- (0,1,1);
      \draw[thick,opacity=1] (0,0,1) -- (1,0,1); 
      \draw[thick,opacity=1] (1,1,0) -- (1,1,1);
      \draw[thick,opacity=1] (0,1,1) -- (1,1,1);
      \draw[thick,opacity=1] (1,0,1) -- (1,1,1);

    \draw[thick,opacity=1] (-0.8,-0.8,-0.8) -- (0,0,0);
    \draw[thick,opacity=1] (2,-0.8,-0.8) -- (1,0,0);
    \draw[thick,opacity=1] (2,2,-0.8) -- (1,1,0);
    \draw[thick,opacity=1] (2,2,2) -- (1,1,1);
    \draw[thick,opacity=1] (2,-0.8,2) -- (1,0,1);
    \draw[thick,opacity=1] (-0.8,2,-0.8) -- (0,1,0);
    \draw[thick,opacity=1] (-0.8,-0.8,2) -- (0,0,1);
    \draw[thick,opacity=1] (-0.8,2,2) -- (0,1,1);

      \draw[dashed,opacity=1] (-0.8,-0.8,-0.8) -- (-0.8,2,-0.8);
      \draw[dashed,opacity=1] (-0.8,-0.8,-0.8) -- (2,-0.8,-0.8);
      \draw[dashed,opacity=1] (-0.8,-0.8,-0.8) -- (-0.8,-0.8,2);
      \draw[thick,opacity=1] (-0.8,-0.8,2) -- (-0.8,2,2);
      \draw[thick,opacity=1] (-0.8,-0.8,2) -- (2,-0.8,2);
      \draw[thick,opacity=1] (-0.8,2,-0.8) -- (-0.8,2,2);
      \draw[thick,opacity=1] (-0.8,2,-0.8) -- (2,2,-0.8);
      \draw[thick,opacity=1] (2,-0.8,-0.8) -- (2,2,-0.8);
      \draw[thick,opacity=1] (2,-0.8,-0.8) -- (2,-0.8,2);
      \draw[thick,opacity=1] (2,-0.8,2) -- (2,2,2);
      \draw[thick,opacity=1] (-0.8,2,2) -- (2,2,2);
      \draw[thick,opacity=1] (2,2,-0.8) -- (2,2,2);

\draw[fill=gray!10] (0,0,0) circle (0.45em); 
\draw[fill=gray!10] (1,0,0) circle (0.45em); 
\draw[fill=gray!10] (0,1,0) circle (0.45em); 
\draw[fill=gray!10] (0,0,1) circle (0.45em); 
\draw[fill=gray!10] (2, -0.8, -0.8) circle (0.45em);

\node at (0,0,0) {$1$};
\node at (0,1,0) {$4$};
\node at (1,0,0) {$3$};
\node at (0,0,1) {$2$};
\node at (2, -0.8, -0.8) {$5$};

      \draw[dashed,opacity=1] (0,-6,-0.5) -- (0,-5,-0.5);
      \draw[dashed,opacity=1] (0,-6,-0.5) -- (1,-6,-0.5);
      \draw[dashed,opacity=1] (0,-6,-0.5) -- (0,-6,0.5);
      \draw[thick,opacity=1] (0,-5,-0.5) -- (0,-5,0.5);
      \draw[thick,opacity=1] (0,-5,-0.5) -- (1,-5,-0.5);
      \draw[thick,opacity=1] (1,-6,-0.5) -- (1,-5,-0.5);
      \draw[thick,opacity=1] (1,-6,-0.5) -- (1,-6,0.5);
      \draw[thick,opacity=1] (0,-6,0.5) -- (0,-5,0.5);
      \draw[thick,opacity=1] (0,-6,0.5) -- (1,-6,0.5); 
      \draw[thick,opacity=1] (1,-5,-0.5) -- (1,-5,0.5);
      \draw[thick,opacity=1] (0,-5,0.5) -- (1,-5,0.5);
      \draw[thick,opacity=1] (1,-6,0.5) -- (1,-5,0.5);

      \draw[fill=gray!10] (0,-6,-0.5) circle (0.45em); 
    \draw[fill=gray!10] (1,-6,-0.5) circle (0.45em); 
    \draw[fill=gray!10] (0,-5,-0.5) circle (0.45em); 
    \draw[fill=gray!10] (0,-6,0.5) circle (0.45em);

    \node at (0,-6,-0.5) {$1$};
    \node at (0,-5,-0.5) {$4$};
    \node at (1,-6,-0.5) {$3$};
    \node at (0,-6,0.5) {$2$};

\draw[right hook->,opacity=1] (0,-4.5,0) -- (0,-2.5,0.15);
\end{tikzpicture}
\end{center}

\end{example}

\begin{lemma}\label{lemma:isolated-lifts}
	Every $n$-lift $C' \in \mathcal{C}[d+n, q_{d,l}(n)]$ of a $k$-isolated configuration $C \in \mathcal{C}[d,l]$ is $(k+n)$-semi-isolated. Furthermore all the configurations in the same connected component as $C'$ in the puzzle graph $puz[d+n, k+n, q_{d,l}(n)]$ are obtained as $n$-lifts of $C$.
\end{lemma}
\begin{proof}
	Fix $k$. Firstly let us prove that the $n$-lifts of $k$-isolated configurations are $(k+n)$-semi-isolated. We proceed by induction on the integer $n \ge 0$. For simplicity, define $(d',l') := (d+n,q_{d,l}(n))$ and $k' := k+n$. Let $C \in \mathcal{C}[d,l]$ be $k$-isolated, then its $0$-lift is $k$-isolated, and we can consider this a degenerate case of $k$-semi-isolated, thus the case $n=0$ trivially holds. Thus we are left to prove that the $(n+1)$-lift of $C$ is $(k'+1)$-semi-isolated, assuming that the $n$-lift is $k'$-semi-isolated. Let us denote as $C_n \in \mathcal{C}[d',l']$ and $C_{n+1} \in \mathcal{C}[d'+1,q_{d',l'}(1)]$ the $n$-lift and $(n+1)$-lift of $C$, respectively and 
	\[
		\iota \co Q^{k'} \rightarrow Q^{k'+1}
	\]
	the embedding of $Q^{k'}$ into $Q^{k'+1}$ sending $C_n$ to $C_{n+1}$. Let $\mathcal A _n$ be the set of stuck tokens of $C_n$, certifying it is a $k'$-semi-isolated configurations. Let $\mathcal B$ be the set of $k'$-dimensional faces $\tau$ of $Q^{d'}$ such that 
	\[
		\tau \cap \mathcal A _n \neq \emptyset
	\]
	For every $\tau \in \mathcal B$, we have $\#(\tau \cap \mathcal A _n) \ge 2$. Indeed, if this was not the case we could slide the single labeled token in $\tau \cap \mathcal A _n$ to another vertex of $\tau$, possibly after a series of $k'$-moves of tokens in  $\mathcal L_{C_n} \setminus \mathcal A _n$, and this would be a contradiction by definition of $\mathcal A_n$.
	Let $\mathcal B'$ be the set of $(k'+1)$-dimensional faces $\sigma$ of $Q^{k'+1}$ containing some $\iota (\tau)$ for $\tau \in \mathcal B$, then the following holds: 
	\begin{equation} \label{count}
	\# (\sigma \cap \iota (\mathcal A _n) ) \ge 2
	\end{equation}
	for every $\sigma \in \mathcal B'$. As the $(k'+1)$-dimensional faces in $\mathcal B'$ are all the faces of that dimension containing elements of $ \iota (\mathcal A _n)$, by Equation \ref{count}, we have that $ \iota (\mathcal A _n)$ certifies that $C_{n+1}$ is $(k'+1)$-semi-isolated.

	The second part of the statement is trivial.
\end{proof}
\section{Unlabeled puzzle graphs}\label{sec:unlabeled puzzle graph}

In this section, we define the unlabeled puzzle graph and prove some of its properties that we will use often in the following. We start by giving its definition. 

\begin{definition}\label{def:unlabeled-configurations}
    Let $\tilde {\mathcal{C}}[d,l]$ be the set of all configurations of $2^d-l$ unlabeled tokens on the vertices of $Q^d$, and let the $(d,k,l)$-th unlabeled puzzle graph, denoted as $\widetilde{puz}[d,k,l]$, the graph whose set of vertices is $\tilde {\mathcal{C}}[d,l]$ and whose set of edges is defined analogously as for the puzzle graph $puz[d,k,l]$. That is
    \[
	\{(C,D) \in \tilde{\mathcal{C}}[d,l] \times \tilde{\mathcal{C}}[d,l] \mid \text {$C$ can be reached from $D$ by a single $k$-move} \}
	\]
    Observe that for every triple $(d,k,l)$ there is a graph homomorphism 
    \[
    u \co puz[d,k,l] \longrightarrow \widetilde{puz}[d,k,l]
    \]
given by sending a labeled configuration to the underlying unlabeled configuration.
\end{definition}

If we only allow slides to adjacent vertices of the $d$-cube, we get different puzzle graphs with fewer edges. We denote these unlabeled and labeled puzzle graphs by $\widetilde{puz}^\prime[d,k,l]$ and ${puz}^\prime[d,k,l]$ respectively. Note that the connectivity does not change, as any move on a $k$-dimensional face can be decomposed into up to $k$ moves to adjacent vertices. We have the following result for the puzzle graphs of labeled and unlabeled configurations with respect to this sliding rule:
    
\begin{lemma}\label{lemma:bipartiteness}
   For $d \ge 3$ and all admissible choices of  $(k,l)$, the puzzle graphs of labeled and unlabeled configurations $\widetilde{puz}^\prime[d,k,l]$ and ${puz}^\prime[d,k,l]$ constructed with respect to the sliding rule above, are bipartite.
\end{lemma}
\begin{proof}
	Consider $Q^d$ as the standard unit cube in $\R^d$ and let $v_i$ for $i \in \{1, \dots , 2^d - l\}$ be the vertices occupied by a token, seen as vectors of $\R^d$. We argue that the function $\varphi \co \tilde{\mathcal{C}}[d,l] \rightarrow \Z / 2 \Z$  defined as:
	\[
	\{v_i\}_{i \in \{1, \dots , 2^d - l\}} \longmapsto \sum_{i=1}^{2^d -l} \sum_{j=1}^d \overline {v_{i,j}}
	\]
	defines a 2-colouring of the vertices of $\widetilde{puz}'[d,k,l]$. Indeed, a $k$-move changes one and only one coordinate of one of the vertices occupied by a token in the configuration, thus it changes the parity of $\sum_{i=1}^{2^d -l} \sum_{j=1}^d v_{i,j}$. Thanks to this observation, the two sets $\varphi^{-1}(\{0\})$ and $\varphi^{-1}(\{1\})$ give a 2-colouring of the vertices of the graph $\widetilde{puz}'[d,k,l]$. 
	
	As any move on a labeled configuration changes the underlying unlabeled configuration, the same proof gives that the labeled puzzle graph $puz'[d,k,l]$ is bipartite.
\end{proof}

Before continuing analyzing the properties of the unlabeled puzzle graphs we introduce a notation that allows to refer to each face of the $d$-cube in a canonical and easy way. We will also use this notation in the next chapter.

\begin{definition} \label{Definition:star-notation}
	Given a vector $v\in \{0,1,\star\}^d$, we define the set 
	\[
	S_v=\{w\in \{0,1\}^d \mid w_i=v_i\text{ for } v_i\in\{0,1\}, 1\leq i\leq d\}.
	\]
	Let $r=\# (\{v_i \mid v_i=\star\})$, then $S_v$ is an $r$-dimensional face of $Q^d \subset \R^d$.
\end{definition}

\begin{definition}
    We say two $(d-1)$-dimensional faces $\sigma$ and $\sigma'$ of a $d$-dimensional cube are parallel if the vectors representing them in the notation introduced in Definition \ref{Definition:star-notation} are $(\star, \dots ,\star, 0, \star, \dots, \star)$ and $(\star, \dots ,\star, 1, \star, \dots, \star)$, where the only non-$\star$ entries are in the same position.
    
    In particular, the set of vertices of $\sigma$ and that of $\sigma'$ do not intersect, and each vertex in the first set is connected to a vertex in the second one, by an edge of $Q^d$.
\end{definition}

\begin{remark}
    In the following we will use the notation $\sigma \cup \sigma'$ for the $d$-dimensional cube which consists of two parallel $d-1$ faces $\sigma$ and $\sigma'$ and all the edges connecting their vertices.
\end{remark}

\begin{lemma} \label{Lemma:vertices}
	Given two adjacent vertices $v$ and $v'$ of $Q^{k+1}$, starting from the configuration $C \in \tilde{\mathcal{C}}[k+1,k,2^{k+1} - (k+2)]$ where the tokens are at $v'$ and at its adjacent vertices, we can slide $k-1$ times to reach the analogue configuration where the tokens are at $v$ and at its adjacent vertices.
\end{lemma}
\begin{proof}
    Without loss of generality we can assume that $v = (0, \dots 0,)$ and $v' = e_1 = ( 1 , 0 , \dots , 0)$.

	We observe, that the labels at $(0,\dots,0)$ and $e_1$ are already in a right position. Because all subcubes $(\star,\dots,\star, 1, \star,\dots, \star)$ with a 1 in the $k$-th entry contain only one token, we can move the token from $e_k$ to $e_1 + e_k$. This way we obtain the desired configuration.
\end{proof}

\begin{corollary}
    By applying multiple times Lemma \ref{Lemma:vertices} above, we get the same result for any two non adjacent vertices $v$ and $v'$.
\end{corollary}

\begin{figure}
\begin{tikzpicture}[scale=1.3, tdplot_main_coords,axis/.style={->},thick]  

      \draw[dashed,opacity=1] (0,0,0) -- (0,1,0);
      \draw[dashed,opacity=1] (0,0,0) -- (1,0,0);
      \draw[dashed,opacity=1] (0,0,0) -- (0,0,1);
      \draw[thick,opacity=1] (0,1,0) -- (0,1,1);
      \draw[thick,opacity=1] (0,1,0) -- (1,1,0);
      \draw[thick,opacity=1] (1,0,0) -- (1,1,0);
      \draw[thick,opacity=1] (1,0,0) -- (1,0,1);
      \draw[thick,opacity=1] (0,0,1) -- (0,1,1);
      \draw[thick,opacity=1] (0,0,1) -- (1,0,1); 
      \draw[thick,opacity=1] (1,1,0) -- (1,1,1);
      \draw[thick,opacity=1] (0,1,1) -- (1,1,1);
      \draw[thick,opacity=1] (1,0,1) -- (1,1,1);

    \draw[thick,opacity=1] (-0.8,-0.8,-0.8) -- (0,0,0);
    \draw[thick,opacity=1] (2,-0.8,-0.8) -- (1,0,0);
    \draw[thick,opacity=1] (2,2,-0.8) -- (1,1,0);
    \draw[thick,opacity=1] (2,2,2) -- (1,1,1);
    \draw[thick,opacity=1] (2,-0.8,2) -- (1,0,1);
    \draw[thick,opacity=1] (-0.8,2,-0.8) -- (0,1,0);
    \draw[thick,opacity=1] (-0.8,-0.8,2) -- (0,0,1);
    \draw[thick,opacity=1] (-0.8,2,2) -- (0,1,1);

      \draw[dashed,opacity=1] (-0.8,-0.8,-0.8) -- (-0.8,2,-0.8);
      \draw[dashed,opacity=1] (-0.8,-0.8,-0.8) -- (2,-0.8,-0.8);
      \draw[dashed,opacity=1] (-0.8,-0.8,-0.8) -- (-0.8,-0.8,2);
      \draw[thick,opacity=1] (-0.8,-0.8,2) -- (-0.8,2,2);
      \draw[thick,opacity=1] (-0.8,-0.8,2) -- (2,-0.8,2);
      \draw[thick,opacity=1] (-0.8,2,-0.8) -- (-0.8,2,2);
      \draw[thick,opacity=1] (-0.8,2,-0.8) -- (2,2,-0.8);
      \draw[thick,opacity=1] (2,-0.8,-0.8) -- (2,2,-0.8);
      \draw[thick,opacity=1] (2,-0.8,-0.8) -- (2,-0.8,2);
      \draw[thick,opacity=1] (2,-0.8,2) -- (2,2,2);
      \draw[thick,opacity=1] (-0.8,2,2) -- (2,2,2);
      \draw[thick,opacity=1] (2,2,-0.8) -- (2,2,2);

\draw[fill=gray!10] (0,0,0) circle (0.45em); 
\draw[fill=gray!10] (1,0,0) circle (0.45em); 
\draw[fill=gray!10] (0,1,0) circle (0.45em); 
\draw[fill=gray!10] (0,0,1) circle (0.45em); 
\draw[fill=gray!10] (-0.8, -0.8, -0.8) circle (0.45em); 

\node at (-0.1,-0.2,0.2) {$v$};
\node at (0.9,-0.2,0.2) {$v'$};

\end{tikzpicture}
    \caption{The unlabeled configuration considered in Lemma \ref{Lemma:vertices} with $k=3$, $v=(0, \dots, 0)$ and $v'=e_1$}
    \label{Figure:standard-conf}
\end{figure}
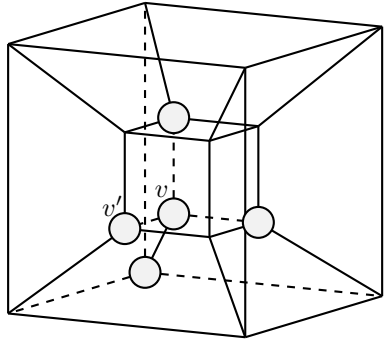

\begin{definition}\label{def:restriction}
    Let $C$ be an unlabeled or labeled configuration on $Q^d$ and $\sigma$ a $ d^\prime$-dimensional face of $Q^d$ with $d^\prime \leq d$. We define the restriction of $C$ to $\sigma$, denoted $C_{\sigma}$, as the configuration on $\sigma$ given by the tokens placed on its vertices.
\end{definition}

We will use this restriction of configurations to specific faces of $Q^d$ with restrictions on the rules for moving the tokens inside of $\sigma$ in the following way: if we are looking at a $k$-rule, tokens on $C_{\sigma}$ can only be moved if they are in a $k$-free state on a $k$-dimensional face of $\sigma$.

\begin{lemma}\label{Lemma:restriction-non-semi-isolated}
    Let $d\geq 3$ and $k=d-1$. Given any configuration $C$ that is $k$-mobile, there is a $k$-dimensional face $\sigma$ with a single token $t$ on it, as in particular $C$ is not isolated. The restriction $C_{\sigma'}$ of $C$ to the $k$-dimensional face $\sigma'$ that is parallel to $\sigma$, is $(k-1)$-mobile.
\end{lemma}

\begin{proof}
First, observe that $\sigma^{\prime }$ is unique, as $k=d-1$. Let us assume by contradiction that $C_{\sigma^{\prime}}$ is $(k-1)$-semi-isolated and denote as $\mathcal A$ the set of vertices certifying it. For any vertex $v$ in $\mathcal A$,  any $(k-1)$-dimensional face of $\sigma^{\prime}$ having $v$ as a vertex, contains at least another token of $\mathcal A$, as no element of the set $\mathcal A$ can be slid.
Any $k$-dimensional face having $v$ as a vertex contains one of the $(k-1)$-dimensional faces above, having at least another vertex in $\mathcal A$. Thus, the token in $v$ can not slide on any $k$-dimensional face, which implies that no token on the vertices of the set $\mathcal A$ can be slid. This, in turn, means that $C$ is $k$-semi-isolated, which contradicts our hypothesis.
\end{proof}

The same method of splitting the $d$-dimensional cube into two parallel $(d-1)$-dimensional faces can be used in the general case to prove the following Lemma.

\begin{lemma}\label{lemma:splitting}
    Let $d\geq 3$, $1 < k < d$, and assume that $l$ is such that there exists a $k$-mobile unlabeled configuration $C$ in $\widetilde{puz}[d,k, l]$. Then, we can choose a $(d-1)$-dimensional face $\sigma$, such that $C_\sigma$ is $k$-mobile and the restriction $C_{\sigma'}$ to the parallel face $\sigma'$ is $(k-1)$-mobile.
\end{lemma}

\begin{figure}
\begin{tikzpicture}[scale=1.3, tdplot_main_coords,axis/.style={->},thick]  

 \draw[dashed,opacity=1] (0,0,0) -- (0,1,0);
      \draw[dashed,opacity=1] (0,0,0) -- (1,0,0);
      \draw[dashed,opacity=1] (0,0,0) -- (0,0,1);
      \draw[thick,opacity=1] (0,1,0) -- (0,1,1);
      \draw[thick,opacity=1] (0,1,0) -- (1,1,0);
      \draw[thick,opacity=1] (1,0,0) -- (1,1,0);
      \draw[thick,opacity=1] (1,0,0) -- (1,0,1);
      \draw[thick,opacity=1] (0,0,1) -- (0,1,1);
      \draw[thick,opacity=1] (0,0,1) -- (1,0,1); 
      \draw[thick,opacity=1] (1,1,0) -- (1,1,1);
      \draw[thick,opacity=1] (0,1,1) -- (1,1,1);
      \draw[thick,opacity=1] (1,0,1) -- (1,1,1);

\draw[fill=gray!10] (0,0,0) circle (0.45em); 
\draw[fill=gray!10] (1,0,0) circle (0.45em); 
\draw[fill=gray!10] (0,0,1) circle (0.45em); 
\draw[fill=gray!10] (0,1,0) circle (0.45em);

\node at (0,1,0) {$3$};
\node at (0,0,0) {$1$};
\node at (1,0,0) {$2$};
\node at (0,0,1) {$4$};

\draw[right hook->,opacity=1] (0,-2.8,0) -- (0,-1.3,0.15);
     
\tikzset{shift={(0,-3.5,-0.3)}}

      \draw[thick,opacity=1] (0,0,0) -- (0,0,1);
      \draw[thick,opacity=1] (0,0,0) -- (1,0,0);
      \draw[thick,opacity=1] (1,0,0) -- (1,0,1);
      \draw[thick,opacity=1] (1,0,1) -- (0,0,1);

\draw[fill=gray!10] (0,0,0) circle (0.45em); 
\draw[fill=gray!10] (1,0,0) circle (0.45em); 
\draw[fill=gray!10] (0,0,1) circle (0.45em); 

    \node at (0,0,0) {$1$};
    \node at (1,0,0) {$2$};
    \node at (0,0,1) {$4$};

\tikzset{shift={(0,8,0.6)}}

    \draw[thick,opacity=1] (0,0,0) -- (0,0,1);
    \draw[thick,opacity=1] (0,0,0) -- (1,0,0);
    \draw[thick,opacity=1] (1,0,0) -- (1,0,1);
    \draw[thick,opacity=1] (0,0,1) -- (1,0,1);
    
    \draw[fill=gray!10] (0,0,0) circle (0.45em); 

    \node at (0,0,0) {$3$};
    
    \draw[left hook->,opacity=1] (0,-1.3,0.275) -- (0,-2.8,0.125);
\end{tikzpicture}
 \caption{The splitting of a 2-mobile configuration on the 3-cube into a 1-mobile (left) and 2-mobile (right) configurations of the 2-cube.}
    \label{Figure:splitting-conf}
\end{figure}
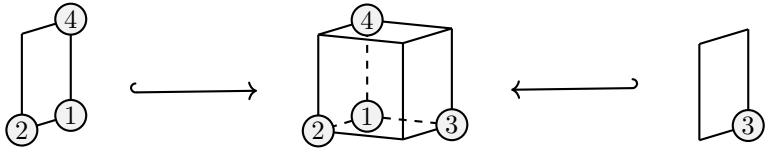
The following proposition follows from the previous Lemmas.

\begin{prop}\label{Proposition:configuration-connected-k}
The subgraph of $k$-mobile configurations of the unlabeled puzzle graph $\widetilde{puz}[k+1, k,  2^{k+1}-(k+2)]$ is connected.
\end{prop}

\begin{proof}
We prove this by induction on $k$. For $k=1$ this trivially holds.

For $k \geq 2$, let $C$ be the configuration with tokens at $v=(0, \dots, 0)$ and its adjacent vertices. For an example of such a configuration $C$, for $k=2$ have a look at Figure \ref{Figure:standard-conf}. 

Let now $D$ be any $k$-mobile configuration, then there is a $k$-dimensional face $\sigma$ with a single token $t$ on it. The restriction $D_{\sigma'}$ to the $k$-dimensional face $\sigma'$, parallel to $\sigma$ is $(k-1)$-mobile thanks to Lemma \ref{Lemma:restriction-non-semi-isolated}, because $D$ is $k$-mobile.
Thanks to Lemma \ref{Lemma:vertices} starting from position $C$, with a series of slides we can reach a configuration where the token $t$ is on the same $k$-dimensional face it occupies in configuration $D$ and all the other tokens are on the parallel $k$-dimensional face $\sigma'$. Let us call this configuration $E$ and its restriction to $\sigma'$ as $E_{\sigma'}$.

We can reproduce the $(k-1)$-rule on $\sigma'$. Indeed, any $(k-1)$-dimensional face $\tau'$ of $\sigma'$ forms a $k$-dimensional face $\rho$ together with the vertices of a $(k-1)$-dimensional face $\tau$ of $\sigma$. 
As $t$ can be moved freely on $\sigma$, we can move it out of $\tau$. Thus if there is only one token on $\tau'$ we can move it on $\rho$. In particular, it can move on $\tau'$. 
By induction we can transform $E_{\sigma'}$ to $D_{\sigma'}$ and thus $E$ to $D$.

Since $C$ is connected to $E$, it follows that $C$ is connected to $D$ in $\widetilde{puz}[k+1,k, 2^{k+1}-(k+2)]$.
\end{proof}

To prove the general case, we use a similar method and prove that the composition of a $(k-1)$-mobile and a $k$-mobile configuration will again be $k$-mobile.
    
\begin{prop}\label{Proposition:composition/unlabeled-connectivity}
For any admissible pair of integers $(d,k)$, the following hold:
 \begin{enumerate}[(i)]
     \item \label{point:composition} let $C$ and $D$ be respectively a $k$-mobile configuration and a $(k-1)$-mobile configuration on $Q^{d-1}$. Every configuration $E$ on $Q^d$, such that there exist two parallel $(d-1)$-dimensional faces $\sigma$ and $\sigma'$ for which $E_{\sigma} = C$ and $E_{\sigma'} = D$, is $k$-mobile;
     \item \label{point:unlabeled-connectivity-general} the subgraph of the unlabeled puzzle graph $\widetilde{puz}[d, k, l]$ whose vertices are all the $k$-mobile configurations is connected.
 \end{enumerate}
\end{prop}

\begin{proof}
    Let $d = k + n$. If $l$ is such that $\widetilde{puz}[d,k,l]$ contains no $k$-mobile configurations, the statement holds by vacuity.
    We prove the statement by induction over $k$ and $n$. 
    Clearly $\widetilde{puz}[d, 1, l]$ is connected for every value of $l$ and the subgraph of $\widetilde{puz}[k, k, l]$ whose vertices are all the $k$-mobile configurations is connected.
    Also the statement holds for $d = k + n < 4$ as shown in the computations 
    Now let us assume (\ref{point:composition}) and (\ref{point:unlabeled-connectivity-general}) hold for all $k' \leq k$, $n' \leq n$ with $k' + n' < k + n$ and let us prove statement (\ref{point:composition}) for $k$ and $n$.

    Let $C$ be a configuration and $\sigma$ and $\sigma'$ parallel $(k+n-1)$-dimensional faces such that $C_{\sigma}$ is $k$-mobile and $C_{\sigma'}$ is $(k-1)$-mobile.
    Let $\mathcal A_C$ be the set of stuck tokens of $C$.
    As $C_{\sigma}$ is $k$-mobile, $\mathcal A_C$ can not intersect the set $\mathcal L_{C_{\sigma}}$ of tokens in $C_{\sigma}$.
    
    In $C$, we can reproduce the $(k-1)$-rule on $\sigma'$. Indeed, let $t$ be a token of $C_{\sigma'}$ that is in a free $(k-1)$-state and let $\tau'$ the $(k-1)$-dimensional face of $\sigma'$ containing no other token but $t$.
    Then there exists a $(k-1)$-dimensional face $\tau$ of $\sigma$ such that $\tau \cup \tau'$ is $k$-dimensional face of $Q^d$.
    As $C_{\sigma}$ is $k$-mobile, by inductive hypothesis on statement (\ref{point:unlabeled-connectivity-general}), we can perform $k$-moves on $\sigma$ to reach from $C$ a configuration $\hat{C}$ whose restriction to $\sigma$ is any $k$-mobile configuration on $\sigma$. In particular we can move in such a way that $\tau$ is empty.
    The face $\tau \cup \tau'$ then shows that $t$ is in a free $k$-state. So $t$ can perform any $k$-move on $\tau \cup \tau'$ including all $(k-1)$-moves on $\tau'$.
    If we do not move $t$ out of $\tau'$, we can perform the moves on $\sigma$ we needed to free up $\tau$ in reversed order. This sequence of moves corresponds to a simple $(k-1)$-move on $C_{\sigma'}$.

    In conclusion, we have that the intersection $\mathcal A_C \cap \mathcal L_{C_{\sigma'}}$ has to be empty, as $C_{\sigma'}$ is $(k-1)$-mobile by hypothesis and we can perform any $(k-1)$-move we could perform on $C_{\sigma'}$, as a $k$-move on $C$.
    Thus $\mathcal A_C$ is empty and $C$ is $k$-mobile. This completes the proof of statement (\ref{point:composition}).
    
    To show statement (\ref{point:unlabeled-connectivity-general}) holds for $k$ and $n$ we have to prove that given two unlabeled $k$-mobile configurations $D$ and $E$, it is possible to reach $E$ from $D$ with a series of slides.
    
    Let us assume that for both $D$ and $E$ we can choose the same parallel $(k+n-1)$-dimensional faces $\sigma$ and $\sigma'$ such that the restrictions $D_{\sigma}$ and $E_{\sigma}$ are $k$-mobile and the restrictions $D_{\sigma'}$ and $E_{\sigma'}$ are $(k-1)$-mobile. (see Lemma \ref{lemma:splitting})

    In the same way as above we are again able to reproduce the $(k-1)$-rule on $\sigma'$.
    As $D_{\sigma'}$ is $(k-1)$-mobile by inductive hypothesis on statement (\ref{point:unlabeled-connectivity-general}), we can transform it into $E_{\sigma'}$ using the $(k-1)$-rule on the $(k+n-1)$-dimensional face $\sigma'$. Analogously, $D_{\sigma}$ is $k$-mobile and we can transform it into $E_{\sigma}$.
    By this we first transformed $D_{\sigma'}$ into $E_{\sigma'}$ and then $D_{\sigma}$ into $E_{\sigma}$.
    Thus in this case $D$ is in the same connected component of $\widetilde{puz}(k+n,k ,l) $ as $E$.
    
    Let now $D$ such that the restriction $D_{\sigma}$ to a $(k+n-1)$-dimensional face $\sigma$ is $k$-mobile and the restriction $D_{\sigma'}$ to its parallel $(k+n-1)$-dimensional face $\sigma'$ is $(k-1)$-mobile.
    Let also $E$ such that the restriction $E_{\rho}$ to a $(k+n-1)$-dimensional face $\rho$ is $k$-mobile and the restriction $E_{\rho'}$ to its parallel $(k+n-1)$-dimensional face $\rho'$ is $(k-1)$-mobile.
    The case $\rho = \sigma$ is already proven, so what is left to show is that $D$ is in the same connected component of $\widetilde{puz}(d,k,l)$ as a configuration $\hat{D}$, such that the restriction $\hat{D}_{\rho}$ is $k$-mobile and the restriction $\hat{D}_{\rho'}$  is $(k-1)$-mobile.

    If $\rho$ is not parallel to $\sigma$, we have: 
    \[
    \rho = (\sigma \cap \rho) \cup (\sigma' \cap \rho)
    \]
    and 
    \[
    \rho' = (\sigma \cap \rho') \cup (\sigma' \cap \rho')
    \]
    where $\sigma \cap \rho$, $\sigma' \cap \rho$, $\sigma \cap \rho'$ and $\sigma' \cap \rho'$ are  $(k+n-2)$-dimensional faces. As we can again reproduce the $(k-1)$-rule on $\sigma'$, thanks to the inductive hypothesis on statement (\ref{point:unlabeled-connectivity-general}), we can perform $(k-1)$-moves on $\sigma'$ in such a way that the restriction to $(\sigma' \cap \rho')$ is $(k-2)$-mobile configuration and the restriction to its parallel $(k-2)$-face $(\sigma' \cap \rho)$ is $(k-1)$-mobile.
    
    Again by inductive hypothesis on Statement (\ref{point:unlabeled-connectivity-general}), we can perform $k$-moves on $\sigma$  in such a way that the restriction to $(\sigma \cap \rho')$ is $(k-1)$-mobile and the restriction to its opposite face $(\sigma \cap \rho)$ is $k$-mobile.
    Thus we reached a configuration $\hat{D}$ such that the restriction $\hat{D}_{\rho}$ is the composition of the $k$-mobile configuration $\hat{D}_{(\sigma \cap \rho)}$ and the $(k-1)$-mobile configuration  $\hat{D}_{(\sigma' \cap \rho)}$ and thus $\hat{D}_{\rho}$ is $k$-mobile by statement  (\ref{point:composition}).
    Moreover the restriction $\hat{D}_{\rho'}$ is the composition of the $(k-1)$-mobile configuration $\hat{D}_{(\sigma \cap \rho')}$ and the $(k-1)$-mobile configuration  $\hat{D}_{(\sigma' \cap \rho')}$ and thus $\hat{D}_{\rho'}$ is $(k-1)$-mobile again by statement
    (\ref{point:composition}).
    Thus $\hat{D}$ has the desired properties.

    In the case that $\rho$ is parallel to $\sigma$, performing the same procedure twice yields the desired result. 
    Thus, any two unlabeled $k$-mobile configurations belong to the same component in $\widetilde{puz}[d, k, l]$ 
\end{proof}

\section{Explicit formulas for $S(d,k)$: proofs of Theorem \ref{theorem:value-of-L(k+1,k)} and \ref{theorem:d=k+n-intro}}\label{sec:values-of-S(d,k)}

In this section, we give a general formula for $S(d,k)$ for any admissible $(d,k)$. We start by introducing Lemma \ref{lemma:no-empty-faces} to prove in Theorem \ref{theorem:equation for L(k+1,k)} that $
S(k+1,k)=2^{k+2}-(k+2)$.
\begin{lemma}\label{lemma:no-empty-faces}
For $k > 1$, every configuration in $puz[k+1, k, 2^{k+1} - (k + 2)]$  that has an empty $k$-face is $k$-isolated or $k$-semi-isolated. 
Furthermore, in $puz[k+1, k, 2^{k+1} - (k + 3)]$ every configuration is $k$-isolated or $k$-semi-isolated.
\end{lemma}
\begin{proof}
We prove the statement by induction.
For $k=2$ both statement hold, as any configuration on the $3$-cube using the $2$-rule with all 4 labels on one $2$-face is isolated. Similarly, any configuration on the $3$-cube using the $2$-rule with 5 labels is either semi-isolated, if 4 labels are stuck on one $2$-face and the last label can move on the parallel face, or isolated.

Let now $k \ge 3$ and assume that both statements are true for $k-1$.
To prove the first statement, let $D$ be a configuration with $2^{k+1} - (k + 2)$ unoccupied vertices and with an empty $k$-dimensional face $\sigma$. If the restriction of $D$ to the parallel $k$-face $\sigma'$ is $(k-1)$-isolated, then $D$ is $k$-isolated. Otherwise there exists a label $x$ that can move on a $(k-1)$-face of $\sigma'$. Because this $(k-1)$-face forms a $k$-face of the whole $(k+1)$-cube together with $2^k-1$ empty vertices of $\sigma$, starting from $D$ we can move $x$ out of $\sigma'$. Thus $D$ can be seen as a $1$-lift of a configuration on $\sigma'$, which is $(k-1)$-semi-isolated by induction, thus by Lemma \ref{lemma:isolated-lifts}, $D$ is $k$-semi-isolated.

The second statement immediately follows from Lemma \ref{lemma:isolated-lifts} as a configuration with $2^{k+1} - (k + 3)$ unoccupied vertices, is either isolated or can be seen as a lift of a configuration on the $k$-cube with $l := 2^{(k-1)+1} - ((k- 1) + 3)$ unoccupied vertices, which is $(k-1)$-semi isolated by induction (notice that in this case $q_{k,l}(1) = 2^{k+1} - (k + 3)$).
\end{proof}

\begin{theorem}\label{theorem:equation for L(k+1,k)}
	For any $k \ge 1$, the following equality holds:
	\begin{equation}
	S(k+1,k) = \sum_{i = 1}^{k} (2^i - 1) = 2^{k+1} - (k + 2)
	\end{equation}
	i.e.\ the even solvability property holds with at most $k+2$ labeled tokens, by definition of $S(k+1,k)$.
\end{theorem}

\begin{proof}
	We prove the statement by induction on the natural number $k$. Computations show that the statement holds for $k = 1,2,3,4$.

	Let $Q^{k+1}$ be embedded in $\R^{k+1}$. Consider the component that contains the configuration $C$ with the $k+2$ labeled tokens at the vertices $(0,\dots,0)$ and $e_1, \dots , e_{k+1}$, where $e_i$ is the $i$-the standard basis vector for $\R^{k+1}$. 
	
    Let us denote by $C_{\sigma'}$ the restriction of $C$ to the $k$-dimensional subcube $\sigma' = (\star,\star,\dots,\star,0)$. Arguing in the same way as in the inductive step of the proof of Proposition \ref{Proposition:configuration-connected-k}, we can reproduce the $(k-1)$-rule on $\sigma'$. As $\sigma'$ contains $k+1$ labeled tokens we can use the inductive hypothesis to get that all configurations on $Q^{k+1}$ obtained by even permutations of the $k+1$ labeled tokens of $C_{\sigma'}$ on $\sigma'$ are in the same connected component as $C$. 
    Choosing a different $k$-dimensional face $(\star,\dots,\star,0,\star,\dots,\star)$ as the face $\sigma'$, the same holds for all configurations obtained from $C$ by permuting all its $k+2$ tokens by even permutations.

    Observe that we do not lose generality by starting from this specific configuration $C$, as by Proposition \ref{Proposition:configuration-connected-k} we can transform any unlabeled $k$-mobile configuration to any other unlabeled $k$-mobile configuration, in particular to $u(C)$. Thus by the conjugation Lemma in \cite{hex} even solvability property follows for $puz[k+1,k,2^{k+1}-(k+2)]$ and by Lemma \ref{lemma:no-empty-faces} we get $S(k+1,k) = 2^{k+1} - (k+2) $.
\end{proof}

 We remind the reader that a permutation group $G$ is transitive on a set $X$ if, given any two elements $a,b \in X$  there exists a permutation $\alpha \in G$ such that $\alpha(a)=b$. Such a group $G$ is said to be primitive if  there is no nontrivial partition of $X$ that is invariant under the group action. Before generalizing the previous result to $S(n+k , k)$ for all $n>1$ we state the following lemma that gives a condition for a permutation group to contain the alternating group $A_n$, i.e.\ the group containing all the even permutations of a set of $n$ elements.

\begin{lemma}\label{lemma:alternating-group}
    Let $G$ be a transitive permutation group on a set of $n$ labeled tokens and suppose that $G$ contains a 3-cycle. If $G$ is primitive then it contains the alternating group $A_n$.
\end{lemma}

\begin{proof}
    A proof can be found, for instance, in \cite[Section 3]{wilson}.
\end{proof}

For general values of $(d,k)$ with $d>k+1$, we give a general formula for $S(d,k)$ in the following Theorem \ref{theorem:d=k+n}. The case $d = k+1$ was analyzed in Theorem \ref{theorem:equation for L(k+1,k)}.
Recall that, as discussed in Section \ref{sec:introduction}, $S(i,i) = 2^{i} - 1$ and $S(i,1) = 1$ for all $i\geq 1$.

\begin{theorem}\label{theorem:d=k+n}
For any $k \ge 2$, $n \ge 1$ we have the following recursive formula:  
\begin{equation}\label{formula:recursive_d=k+n}
S(n+k , k)  = S(n+k - 1, k-1) + S(n+k - 1, k) 
\end{equation}

Furthermore, the following statements hold:
    \begin{enumerate}[(i)]

        \item \label{one} every configuration in $puz[n+k, k, S(n+k,k) - 1]$ is $k$-semi-isolated or $k$-isolated;  
        
        
        \item \label{three} the graph $puz[n+k,k, S(n+k,k)]$ has the even solvability property.
    \end{enumerate}

\end{theorem}

\begin{proof}
We prove together statements (\ref{one}) and (\ref{three}) and Formula (\ref{formula:recursive_d=k+n}) by an inductive argument on $k$ and $n$.

For the base case of Formula (\ref{formula:recursive_d=k+n}) we know that $S(n+1,1) = 1$ holds by \cite{wilson} and that $S(k+1,k) = 2^{k+1} - k - 2$ by Theorem \ref{theorem:equation for L(k+1,k)}. The base case of (\ref{one}) is given by Lemma \ref{lemma:no-empty-faces} and the base case for (\ref{three}) is Theorem \ref{theorem:equation for L(k+1,k)}.

We now assume that Formula (\ref{formula:recursive_d=k+n}) and both statements (\ref{one}) and (\ref{three}) hold for all $k' \leq k$, $n' \leq n$ with $k' + n' < k + n$, and let us proceed to prove they hold for $k$ and $n$.

To prove statement (\ref{one}) we proceed by contradiction. Let us assume there exists a $k$-mobile configuration $D$ in $puz[k+n, k, S(k+n,k) - 1]$. By Lemma \ref{lemma:splitting} there exist a $(k+n-1)$-dimensional face $\sigma$ such that the restriction $D_{\sigma}$ is $k$-mobile. By inductive hypothesis, $D_{\sigma}$ will have at least $S(k+n-1, k)$ unoccupied vertices. Thus the $(k+n-1)$-dimensional face $\sigma'$  parallel to $\sigma$ will have less than $S(k+n-1, k-1)$ unoccupied vertices and thus the restriction $D_{\sigma'}$ is $(k-1)$-semi-isolated. This contradicts  the inductive hypothesis on (\ref{one}). This proves that (\ref{one}) holds.

To prove statement (\ref{three}), we will show that we can apply any even permutation to a $k$-mobile configuration by a series of $k$-moves.

Let $D$ be a $k$-mobile configuration in $puz[k+n,k,S(k+n,k)]$ and $\sigma$ a $(k+n-1)$-dimensional faces such that the restriction to $D_{\sigma}$ is $k$-mobile.
Note that by inductive hypothesis, statement (\ref{three}) holds for the restriction to $D_{\sigma}$. Thus we can permute the labels on $\sigma$ by any even permutation.
By Proposition \ref{Proposition:composition/unlabeled-connectivity}, statement (\ref{point:unlabeled-connectivity-general}) we can empty (by consecutive $k$-moves) any $(k-1)$-dimensional face of $\sigma$. This allows us to perform any $(k-1)$-move of any token on the face $\sigma'$, that is, the face parallel to $\sigma$, as described in the aforementioned proof of Proposition \ref{Proposition:composition/unlabeled-connectivity}. 

Using this, the inductive hypothesis on the  $(k-1+n)$-dimensional face $\sigma'$, and noticing that $\sigma'$ has $S(k+n-1, k-1)$ unoccupied vertices, we can obtain all labeled configurations on $\sigma'$ by even permutations of the labels in $D_{\sigma'}$.
In other words, what is left to prove is that the permutation group generated by slides is transitive, i.e.\ that we can move tokens between $\sigma$ and $\sigma'$.

As $D_{\sigma'}$ is $k$-mobile, by Lemma \ref{lemma:splitting} there exists a $(k+n-2)$-dimensional face $\tau'$ of $\sigma'$ such that the restriction $D_{\tau'}$ is $(k-1)$-mobile. 
There exist also a $(k+n-2)$-dimensional face $\tau$ of $\sigma$, such that $\tau \cup \tau'$ is a $(k+n-1)$-dimensional face $\rho$.
 
We want to use $\rho$ to move tokens between $\sigma$ and $\sigma'$. By the conjugation Lemma of \cite{hex} what we need to prove is that by performing slides on $\sigma$ from $D$ we can reach a configuration $E$ such that: 1)  $E_\rho$ is $k$-mobile, and 2)  $E_\rho$ contains $S(k+n-1,k)$ unoccupied vertices. 

We first prove 1). Observe that $D_\sigma$ is $k$-mobile. Thus, by point (\ref{point:unlabeled-connectivity-general}) of Proposition \ref{Proposition:composition/unlabeled-connectivity} we can perform slides on $\sigma$ to reach a configuration $E$ such that $E_\tau$ is $k$-mobile.
Then, by point (\ref{point:composition}) of the same proposition, as $\rho = \tau \cup \tau'$, the restriction $E_\rho$ is $k$-mobile. 

Then, the restriction $E_\tau$ has $S(k+n-2,k)$ unoccupied vertices and the restriction $E_{\tau'}$ has $S(k+n-2,k-1)$ many. Thus, we can use the inductive hypothesis on Formula (\ref{formula:recursive_d=k+n}) to get that $E_\rho$ has $S(k+n-1, k)$ unoccupied vertices, proving 2). 

Finally, by using the inductive hypothesis in the statement (\ref{three}) we conclude that $E_\rho$ can be permuted by any even permutation. This allows us to move any token from $\sigma'$ to $\sigma$ and viceversa, as we wanted. 

Thus we can perform any even permutation on $D_\sigma'$ and on $D_\sigma$ and
move tokens between $\sigma'$ and $\sigma$ via $k$-moves, so the group of permutations on the labeled tokens is transitive and primitive.
As it also contains a $3$-cycle, by Lemma \ref{lemma:alternating-group} we can apply any even permutation to $D$, proving statement (\ref{three}). For $S(n+k,k)$ formula (\ref{formula:recursive_d=k+n}) holds.
\end{proof}

\begin{remark}
Notice that, for $n=1$, we recover the formula of Theorem \ref{theorem:equation for L(k+1,k)}.
\end{remark}

\begin{corollary}
    The following equality holds, for any $k \ge 2$ and $n \ge 2$:
    \begin{equation}\label{formula:d=k+n}
    S(n+k , k) = \sum_{i = 2}^{k} \binom{n +k-i}{k-i} (2^{i} - 1) + \sum_{i = 2}^{n} \binom{n+k-i}{k-1}. 
    \end{equation}
\end{corollary}

\begin{proof}
    We apply iteratively formula (\ref{formula:recursive_d=k+n}) until only terms of the form $S(i,i)$ and $S(i,1)$ for $i \geq 2$ appear. We need to understand which terms appear and their coefficients.
    Note that at each iteration the first argument of $S(\cdot, \cdot)$ reduces by 1 and the second argument stays the same in one summand and reduces by 1 in the other.
    For a term $S(i,i)$ to appear in this way, it is necessary to iterate the formula $(n + k - i)$ times. The right hand side reaches $i$ in exactly $k-i$ of the branches. Thus the term $S(i,i)$ arises $\binom{n+k-i}{k-i}$ times.
    On the other hand, for the term $S(i,1)$ to appear it is necessary to iterate the formula $(n + k - i)$ times. The right hand side decreases to $1$ in exactly $k - 1$ of the branches. Thus the term $S(i,1)$ arises $\binom{n+k-i}{k - 1}$ times. 
    Since $S(i,i) = 2^i -1$ and $S(i,1) = 1$ for all $i$, the statement follows.
\end{proof}

\section{$k$-mobile configurations and structure of the puzzle graph}\label{sec:k-mobile-configurations}
We state again Theorem \ref{thm:S(d,k) = smallest-containing-not-k-semi-isolated} and Theorem \ref{theorem:maximality-of-L-introduction} and proceed with their proofs. 

\begin{theorem}
 For any admissible tuple of integers $(d, k)$, the value $S(d,k)$ is the smallest number of unoccupied vertices such that $puz[d, k, S(d,k)]$ contains a $k$-mobile configuration. 
 Moreover, the puzzle graph $puz[ d , k , S(d,k) ]$ has either one or two connected components containing $k$-mobile configurations.     
\end{theorem}

\begin{proof}
The first part of the result follows from Lemma \ref{lemma:no-empty-faces} and statement (\ref{one}) of Theorem \ref{theorem:d=k+n} above. As it is shown that for every $l < S(d,k)$ all configurations are $k$-isolated or $k$-semi-isolated. From Lemma \ref{lemma:even-solvability-k-mobile} it follows that $puz[d,k,S(d,k)]$ contains $k$-mobile configurations. Thus $S(d,k)$ is the smallest number of unoccupied vertices that allows for $k$-mobile configurations. I.e.\ given any $k$-mobile configuration it is possible to obtain any other configuration that arises as an even permutation on its labels.

To prove the second part note that the unlabeled puzzle graph $\widetilde{puz}[d,k,S(d,k)]$ is connected as shown in Proposition \ref{Proposition:configuration-connected-k} and statement (\ref{point:unlabeled-connectivity-general}) of Proposition \ref{Proposition:composition/unlabeled-connectivity}.
Thus given any two $k$-mobile labeled configurations $C$ and $D$  their underlying unlabeled configurations $u(C)$ and $u(D)$ are connected, i.e.\ there exists a series of slides transforming $u(C)$ to $u(D)$. 
Applying the same slides to the labeled configuration $C$ leads to a configuration $E$ such that $u(E) = u(D)$. 
If $E$ arises from an even permutation of the labels in $D$, by even solvability property $E$ is connected to $D$ in the puzzle graph, thus $C$ is to $D$.
If $E$ arises from an odd permutation of the labels in $D$, it is connected to $D$ if and only if there exists a configuration that can be permuted in an odd way. This follows directly from the conjugation Lemma in \cite{hex}.
\end{proof}
This leads to the following definition.
\begin{definition}
    If a puzzle graph contains exactly 2 components containing $k$-mobile configurations we say this puzzle graph has the \textit{strong parity property}.
\end{definition}\label{definition:strong parity property}
The concept of strong parity property was introduced in \cite{hex} where it was applied to $2$-dimensional sliding puzzles with hexagonal tiles.

We proceed now with the proof of Theorem \ref{theorem:maximality-of-L-introduction}, first stated in Section \ref{sec:introduction}.

\begin{theorem}\label{theorem:maximality-of-L}
    For any admissible pair of integers $(d, k)$ and $l > S(d,k)$, the subgraph of the puzzle graph $puz[d, k, l]$ whose vertices are all the $k$-mobile configurations is connected.
\end{theorem}

\begin{proof}
    Let $C$ be $k$-mobile in $\mathcal{C}_\text{lab}[d,S(d,k)+1]$.  As the even solvability property holds for $S(d,k)$ and removing one more token can only allow more moves, every configuration obtained by applying a 3-cycle on the labels of $C$ is in the same connected component as $C$ in $puz[d,k,S(d,k)+1]$. We just need to prove that we can find a configuration in the same connected component of $C$ that is obtained by applying a transposition on the labels of $C$.
	
	As $C$ is $k$-mobile, we can assume at least two token can slide, indeed, if this was not the case and only one token $t$ could slide, then after a series of slides, at least two will be allowed to slide, as otherwise we would have $C$ is $k$-semi-isolated with $\mathcal A _C = \mathcal L _C \setminus \{t\}$, which is a contradiction. 
	Furthermore, we can assume there are at least two tokens on adjacent $k$-faces of $Q^d$ after a series of slides: assume by contradiction this would not be the case, so that all the tokens that can slide are on non-adjacent $k$-faces after any series of slides. Choose two of these tokens, say $t$ and $t'$, then any configuration obtained from $C$ by a 3-cycle sending $t$ to $t'$ could not be obtained, as in order to move $t$ to the position of $t'$ there has to be a token that can be slid on a $k$-face adjacent to that of $t$. But this contradicts the assumption that the even solvability property holds for $C$.
	
	Let $C'$ in $\mathcal{C}_\text{lab}[d,S(d,k)]$ such that by removing a token $c$ from $C'$ we obtain $C$ and choose two other tokens $a, b$ in $C'$. From the previous remarks, we can assume that $a$ and $c$ can slide and lie on adjacent $k$-faces. By applying the 3-cycle $(abc)$ we can slide $a$ in the position previously occupied by $b$. Removing $c$, we are left to prove that we can slide $b$ from the position occupied by $c$ in $C'$ to the position previously occupied by $a$. But this is obvious as $a$ and $c$ occupied vertices of adjacent $k$-faces.

    Because by point (\ref{point:unlabeled-connectivity-general}) of Proposition \ref{Proposition:composition/unlabeled-connectivity} we have that the unlabeled puzzle graph $\widetilde{puz}[d,k,S(d,k)+1]$ is connected, the labeled puzzle graph is also connected, as every permutation of every $k$-mobile configuration is in the same connected component.

    For $l > S(d,k) + 1$ the statement follows as removing even more tokens just allows for more and more movement.
    
\end{proof}

\begin{remark}
It follows from Theorem \ref{theorem:maximality-of-L} above that $S(d,k)$ is the only candidate for $puz[d,k,S(d,k)]$ to have strong parity property.
In this case there are only two connected components containing $k$-mobile configurations, that are distinguished by parity.
\end{remark}

\section{Computational results}\label{sec:computational-results}

In this Section, we summarize the computational results that we obtained for the puzzle graphs $puz[d, k, l]$ for $d = 3, 4$, and $5$. By doing so, we prove Theorems \ref{small-dimension-3} \ref{small-dimension-4} and \ref{small-dimension-4-2}.

We have implemented the algorithms using Python and Julia. The code and the complete output files from our computations can be found at our Github repository \cite{githubcomputations}. The main ingredient in our algorithms is the breadth-first-search (BFS) algorithm. Starting from a configuration, we compute all configurations that can be reached by performing a single \textit{$k$-move}, that is, exactly those that have distance $1$ in the puzzle graph from the starting configuration. Next, in an iterative process, we apply the same procedure to all configurations found in the previous step until no new configurations are found. It is important to check that each newly found configuration was not already added to the list in a previous step. 

In our computations, we used the fact that two different configurations are elements of isomorphic components of the puzzle graph if one can be obtained from the other by applying either a permutation of the labeled tokens or a transformation of the $d$-cube via an element of its symmetric group. Also, if a $k$-mobile configuration exists, we focus on the components containing $k$-mobile configurations and neglect the smaller components containing only isolated or semi-isolated configurations.  

We focus on four properties of the puzzle graphs that we briefly introduce in what follows.

\begin{itemize}
    \item \textbf{Number of configurations in the biggest component.} The biggest connected component is the component with the largest amount of vertices.
    \item \textbf{Diameter.} The diameter of a connected graph, also called God's number, is defined as the maximum length of shortest paths between all the vertices of the graph. Here, we define the diameter $diam$ of a puzzle graph as the diameter of the biggest connected component. Thus, the diameter tells us that starting from any configuration it is possible to reach any other configuration contained in the same component by performing at most $diam$ number of \textit{$k$-moves}. Unsurprisingly, for $d=4$ and $k=2$, due to the humongous amount of vertices in the biggest components of the puzzle graph, we were not able to compute the exact diameter for $5\leq l \leq8$. However, we were able to deduce some bounds using the unlabeled puzzle graphs. It is well known that finding God's number of sequential puzzles such as the 15 Puzzle is highly nontrivial (see for example \cite{ratner1990n2, demaine2018simple}). Nevertheless, a lot of research and computational resources have been put into God's number service \cite{ratner1986finding, rokicki2014diameter, brungger1999parallel, parberry1995real, alpert2020discrete}. 
    \item \textbf{Connectivity regime.} In the case that there are only isolated components we say that the regime is $k$-isolated. If there are only $k$-isolated and at least one $k$-semi-isolated component, we say that the regime is $k$-semi-isolated. We proved in Theorem \ref{thm:S(d,k) = smallest-containing-not-k-semi-isolated} that any puzzle graph can have at most two $k$-mobile components. If there are two $k$-mobile components we say that the strong parity property holds; otherwise, we say that the puzzle graph is connected (even though there might be $k$-isolated or $k$-semi-isolated components).
    \item \textbf{Number of Components.}
    We count the number of components that contain the maximal number of configurations. We know that all these components will be isomorphic.
\end{itemize}

The following three tables show the results of computations for $d = 2, 3$, and $4$ for the case $k = 1$. As discussed in the introduction and \cite{wilson}, in this case, strong parity property holds exactly for $l=1$ and we have connectivity for $l > 1$.

\begin{center}
\begin{tabular}{|c | c | c | c | c |} 
\hline\hline 
$l$ & \makecell{\#configurations in\\biggest component} & $diam$ & connectivity regime & \#components \\ [0.5ex] 
\hline 
1 &	12 & 6 & strong parity & 2	\\ 	 	 	 	2 &	12 & 4 & connected & 1  	\\ 		 	 	 3 &	4  & 2 & connected & 1 	 \\
\hline\hline  
\end{tabular}
\captionof{table}{$d=2$, $k=1$}
\end{center}

\begin{center}
\begin{tabular}{|c | c | c | c | c |} 
	
\hline\hline 
$l$ & \makecell{\#configurations in\\biggest component} & $diam$ & connectivity regime & \#components \\ [0.5ex] 
\hline 
1 &	$\binom{8}{7}\cdot \frac{7!}{2}$ & 19 & strong parity & 2	\\
2 &	$\binom{8}{6}\cdot 6!$ & 18 & connected & 1 \\
3 &	$\binom{8}{5}\cdot 5!$  & 15 & connected & 1 \\ 
4 &	$\binom{8}{4}\cdot 4!$  & 12 & connected & 1 	 	\\ 
5 &	$\binom{8}{3}\cdot 3!$  & 9 & connected & 1 	 	\\ 
6 &	$\binom{8}{2}\cdot 2!$  & 6 & connected & 1 	 	\\ 
7 &	8  & 3 & connected & 1 	 	\\ 
\hline\hline 

\end{tabular}
\captionof{table}{$d=3$, $k=1$\label{table:d3k1}}
\end{center}
\vspace{0.5cm}

\begin{center}
\begin{tabular}{|c | c | c | c | c |} 
\hline\hline 
$l$ & \makecell{\#configurations in\\biggest component} & $diam$ & connectivity regime & \#components \\ [0.5ex] 
\hline 
1 &	\textcolor{red}{$\binom{16}{15}\cdot \frac{15!}{2}$} & \textcolor{red}{unknown} & strong parity & $2$	\\
2 &	\textcolor{red}{$\binom{16}{14}\cdot 14!$} & \textcolor{red}{unknown} & connected & $1$ \\
$\vdots$ &	$\vdots$ & $\vdots$ & $\vdots$ & $\vdots$\\ 
7 & \textcolor{red}{$\binom{16}{9}\cdot 9!$} & \textcolor{red}{unknown} & connected & $1$\\
8 &	518918400 & 32 & connected & $1$\\
9 & 57657600 & 28 & connected & $1$\\
10 & 5765760 & 24 & connected & $1$ \\
11 & 524160 & 20 & connected & $1$\\
12 & 43680 & 16 & connected & $1$\\
13 & 3360 & 12 & connected & $1$\\
14 & 240 & 8 & connected & $1$\\
15 & 16 & 4 & connected & $1$\\

\hline\hline 

\end{tabular}
\captionof{table}{$d=4$, $k=1$\label{table:d4k1}}
\end{center}

We close the discussion for the case $k=1$ with a conjecture on the diameter of these puzzle graphs for when $l \geq 2$. It is well known that computing the diameter of sliding puzzles with the $1$-rule on general graphs is NP-hard \cite{goldreich}. However, based on our computations the following pattern seems to hold:
\begin{conj}
    For any $l \in \{2,3,\dots,2^d\}$ the diameter of $puz[d,1,l]$ is given by $diam=d \cdot (2^d-l)$.
\end{conj}

For $d=3, k=2$ and $d=4, k=3$ we were able to calculate the puzzle graphs for all possible values of $l$. This proves Theorems \ref{small-dimension-3} and \ref{small-dimension-4} stated in the introduction. We observe that in these cases all the different connectivity regimes are achieved as $l$ varies.

\begin{center}
\begin{tabular}{|c | c | c | c | c |} 
	
\hline\hline 
$l$ & \makecell{\#configurations in\\biggest component} & $diam$ & connectivity regime & \#components \\ [0.5ex] 
\hline 
1 &	1 & 0 & $2$-isolated & $\binom{8}{7}\cdot 7!$	\\
2 &	1 & 0 & $2$-isolated & $\binom{8}{6}\cdot 6!$ \\
3 &	4 & 2 & $2$-semi-isolated & $5! \cdot 6$\\ 
4 &	672 & 10 & strong parity & 2 	 	\\ 
5 &	$\binom{8}{3}\cdot 3!$  & 6 & connected & 1 	 	\\ 
6 &	$\binom{8}{2}\cdot 2!$  & 4 & connected & 1 	 	\\ 
7 &	8  & 3 & connected & 1 	 	\\ 
\hline\hline 

\end{tabular}
\captionof{table}{$d=3$, $k=2$\label{table:d3k2}}
\end{center}

\begin{center}
\begin{tabular}{|c | c | c | c | c |} 
\hline\hline 
$l$ & \makecell{\#configurations in\\biggest component} & $diam$ & connectivity regime & \#components \\ [0.5ex] 
\hline 
1 &	1 & 0 & $3$-isolated & $\binom{16}{15}\cdot 15!$	\\
2 &	1 & 0 & $3$-isolated & $\binom{16}{14}\cdot 14!$ \\
3 &	1 & 0 & $3$-isolated & $\binom{16}{13}\cdot 14!$\\ 
4 &	1 & 0 & $3$-isolated & $\binom{16}{12}\cdot 14!$\\
5 &	1 & 0 & $3$-isolated & $\binom{16}{11}\cdot 14!$\\
6 &	1 & 0 & $3$-isolated & $\binom{16}{10}\cdot 14!$\\ 
7 &	8 & 1 & $3$-semi-isolated & $9! \cdot 8$\\
8 &	8 & 1 & $3$-semi-isolated & $ 8 \cdot \binom{8}{7} \cdot  8!$\\
9 &	8 & 1 & $3$-semi-isolated & $8 \cdot \binom{8}{6} \cdot  7! $\\
10 & 48 & 3 & $3$-semi-isolated & $24 \cdot \binom{6}{4} \cdot 4! \cdot 4$ \\
11 & 137280 & 17 & strong parity & $2$\\
12 & 40320 & 10 & connected & $1$\\
13 & 3360 & 6 & connected & $1$\\
14 & 240 & 4 & connected & $1$\\
15 & 16 & 4 & connected & $1$\\

\hline\hline 

\end{tabular}
\captionof{table}{$d=4$, $k=3$\label{table:d4k3}}
\end{center}

In these two cases, Tables \ref{table:d3k2} and \ref{table:d4k3}, we also found that for every starting configuration the final depth $f$ reached by the BFS algorithm turned out the be very close to the diameter of the puzzle graph. In fact, we observed that
$diam \leq f + 3$. 
We expect this to hold true for bigger puzzle graphs because of their inherent symmetry. This fact is indeed remarkable, as in principle $f \leq diam \leq 2 f$.

Another important observation for these two cases is that exactly at $S(d,k)$ unoccupied vertices, strong parity property holds. This happens at $l=4$ for $d=3,k=2$, and at $l=11$ for $d=4,k=3$, respectively.

\begin{center}
\begin{tabular}{|c | c | c | c | c | } 
\hline\hline 
$l$ & \makecell{\#configurations in\\biggest component} & $diam$ & connectivity regime & \#components \\ [0.5ex] 
\hline 
1 &	1 & 0 & $2$-isolated & $\binom{16}{15}\cdot 15!$	\\
2 &	1 & 0 & $2$-isolated & $\binom{16}{14}\cdot 14!$ \\
3 &	4 & 2 & $2$-semi-isolated & $13! \cdot 24$\\ 
4 &	672 & 10 & $2$-semi-isolated & $ \binom{16}{4} \cdot 8! \cdot 16$ 	 	\\ 
5 &	\textcolor{red}{45664819200} & \textcolor{red}{unknown} & \textcolor{red}{strong parity} & \textcolor{red}{ 2} 	 \\ 
6 &	\textcolor{red}{25430630400} & \textcolor{red}{unknown} & \textcolor{red}{connected} & \textcolor{red}{1}  \\ 
7 &	4064256000 & $21\leq diam \leq 42$& connected &  1	\\
8 &	 517708800 & $17\leq diam\leq 34$& connected &  1 \\
9 &	 57657600 & 17 & connected &  1	\\
10 & 5765760  & 14 & connected &  1	\\
11 & 524160 & 11 & connected &  	1 	\\
12 & 43680  & 9 & connected &  	 1	\\
13 & 3360  & 6 & connected &   1	\\
14 & 240 & 4 & connected & 1  	 	\\
15 & 16 & 4 & connected & 1 	 	\\
\hline\hline 

\end{tabular}
\captionof{table}{$d=4$, $k=2$}\label{table:d4k2}
\end{center}

Finally, with the computations reported in Table \ref{table:d4k2}, we prove Theorem \ref{small-dimension-4-2} stated in the Introduction.
Entries of rows 5 and 6 in this table are computed by running computations on the unlabeled puzzle graph and combining these computations with Lemmata \ref{lemma:cycles-1} and \ref{lemma:cycles-2}. 

Before stating the Lemmas, we observe that all permutations of a configuration $C$ arise from looking at the underlying unlabeled configuration $u(C)$, where $u$ is the homeomorphism defined in Definition \ref{def:unlabeled-configurations}, in the sense that each edge in the unlabeled puzzle graph corresponds to a transposition of a token and an unoccupied vertex.

\begin{lemma}\label{lemma:cycles-1}
    If all simple cycles in the unlabeled puzzle graph correspond to even permutations, then the labeled puzzle graph has the weak parity property.
\end{lemma}

\begin{lemma}\label{lemma:cycles-2}
    If all cycles in a cycle base for $\widetilde{puz}(d,k,l)$ given by a spanning tree correspond to even permutations, then $puz[d,k,l]$ has the weak parity property.
\end{lemma}

The proofs of Lemmatas \ref{lemma:cycles-1} and \ref{lemma:cycles-2} follow easily from the definition and basic properties of unlabeled puzzle graphs given in Section \ref{sec:unlabeled puzzle graph}. Checking that these conditions hold on the unlabeled puzzle graph is way faster than computing the puzzle graph in the labeled case. Indeed, computations by this method allowed us to check that $puz[4,2,S(4,2)]$, $puz[5,4,S(5,4)]$, and $puz[5,2,S(5,2)]$ have all the strong parity property.

\begin{remark}
    Not all $k$-semi-isolated components have to be isomorphic.
    For $d=4$, $k=2$, $l=4$ there are 2-semi-isolated components that allow movement on a 3-dimensional face or a 2-dimensional face depending on the placement of the unoccupied vertices. Recall that there are also isolated components as discussed above.
\end{remark}


\bibliography{bibliography}

\begin{thebibliography}{10}

\bibitem{alpert2020discrete}
Hannah Alpert.
\newblock Discrete configuration spaces of squares and hexagons.
\newblock {\em Journal of Applied and Computational Topology}, 4(2):263--280,
  2020.

\bibitem{videogame}
Moritz Beyer, Stefano Mereta, Érika Rold\'an, and Peter Voran.
\newblock Cubix video game, 2023.
\newblock https://www.erikaroldan.net/cubical-sliding-puzzles.

\bibitem{githubcomputations}
Moritz Beyer, Stefano Mereta, Érika Rold\'an, and Peter Voran.
\newblock Github repository, December 2022.
\newblock https://github.com/p-voran/CubicalSlidingPuzzles.

\bibitem{brungger1999parallel}
Adrian Br{\"u}ngger, Ambros Marzetta, Komei Fukuda, and J{\"u}rg Nievergelt.
\newblock The parallel search bench {ZRAM} and its applications.
\newblock {\em Annals of Operations Research}, 90(0):45--63, 1999.

\bibitem{demaine2018simple}
Erik~D Demaine and Mikhail Rudoy.
\newblock A simple proof that the $(n^2- 1)$-puzzle is hard.
\newblock {\em Theoretical Computer Science}, 732:80--84, 2018.

\bibitem{goldreich}
Oded Goldreich.
\newblock Finding the shortest move-sequence in the graph-generalized 15-puzzle
  is {NP}-hard.
\newblock In {\em Studies in complexity and cryptography. Miscellanea on the
  interplay between randomness and computation}, pages 1--5. Springer, 2011.

\bibitem{hordern}
Edward Hordern.
\newblock {\em Sliding Piece Puzzles (Recreations in Mathematics, 4)}.
\newblock Oxford University Press, 1987.

\bibitem{johnson1879notes}
Wm~Woolsey Johnson.
\newblock Notes on the “15” puzzle, part i.
\newblock {\em American Journal of Mathematics}, 2(4):397--399, 1879.

\bibitem{hex}
Ray Karpman and Érika Rold\'an.
\newblock Parity property of hexagonal sliding puzzles.
\newblock {\em arXiv preprint arXiv:2201.00919}, 2022.

\bibitem{parberry1995real}
Ian Parberry.
\newblock A real-time algorithm for the $(n^2- 1)$-puzzle.
\newblock {\em Information Processing Letters}, 56(1):23--28, 1995.

\bibitem{ratner1986finding}
Daniel Ratner and Manfred Warmuth.
\newblock Finding a shortest solution for the $n x n$ extension of the
  15-{P}uzzle is intractable.
\newblock {\em Proceedings of the AAAI Conference on Artificial Intelligence},
  5:168, 1986.

\bibitem{ratner1990n2}
Daniel Ratner and Manfred Warmuth.
\newblock The $(n^2- 1)$-puzzle and related relocation problems.
\newblock {\em Journal of Symbolic Computation}, 10(2):111--137, 1990.

\bibitem{rokicki2014diameter}
Tomas Rokicki, Herbert Kociemba, Morley Davidson, and John Dethridge.
\newblock The diameter of the {R}ubik's {C}ube group is twenty.
\newblock {\em {SIAM} {R}eview}, 56(4):645--670, 2014.

\bibitem{story1879notes}
William~Edward Story.
\newblock Notes on the “15” puzzle, part ii.
\newblock {\em American Journal of Mathematics}, 2(4):399--404, 1879.

\bibitem{wilson}
Richard~M Wilson.
\newblock Graph puzzles, homotopy, and the alternating group.
\newblock {\em Journal of Combinatorial Theory, Series B}, 16(1):86--96, 1974.

\end{thebibliography}
\bibliographystyle{plain}

\vspace{2cm}
\end{document}